\documentclass[final,12pt]{elsarticle}

\usepackage[hyphens]{url}

\usepackage{booktabs,mathrsfs}       
\usepackage{amsfonts}       
\usepackage{nicefrac}       
\usepackage{microtype}      
\usepackage{lipsum}
\usepackage{fancyhdr}       
\usepackage{graphicx}       
\usepackage{amsmath}
\usepackage{amsthm}
\usepackage{amssymb}
\usepackage[dvipsnames]{xcolor}
\usepackage{multicol}
\usepackage[most]{tcolorbox}
\usepackage[textsize=scriptsize]{todonotes}
\graphicspath{{media/}}     
\usepackage{cleveref}
\usepackage{bm}
\usepackage[displaymath]{lineno}
\usepackage{caption,subcaption}

\newtheorem{lemma}{Lemma}[section]
\newtheorem{proposition}{Proposition}
\crefname{lemma}{Lemma}{Lemmas}
\crefname{theorem}{Theorem}{Theorems}
\crefname{equation}{Equation}{Equations}
\crefname{section}{Section}{Sections}
\crefname{subsection}{Section}{Sections}
\crefname{figure}{Figure}{Figures}
\crefname{align}{Equation}{Equations}
\crefname{gather}{Equation}{Equations}
\crefname{proposition}{Proposition}{Propositions}

\newcommand{\mathbbm}[1]{\text{\usefont{U}{bbm}{m}{n}#1}} 

\DeclareMathOperator*{\argmax}{arg\,max}

\journal{Signal Processing Journal}

\begin{document}

\begin{frontmatter}


\title{Performance Bounds for LASSO under Multiplicative LogNormal Noise: Applications to Pooled RT-PCR Testing}
\author[label1]{Richeek Das}
\ead{richeek@cse.iitb.ac.in}
\affiliation[label1]{organization={Indian Institute of Technology Bombay},
            addressline={Department of Computer Science and Engineering}}
\affiliation[label3]{organization={Indian Institute of Technology Bombay},addressline={Department of Electrical Engineering}}


\author[label3]{Aaron Jerry Ninan}
\ead{aaronjerry12@gmail.com}

\author[label1]{Adithya Bhaskar}
\ead{adithyabhaskar@cse.iitb.ac.in}

\author[label1]{Ajit Rajwade}
\ead{ajitvr@cse.iitb.ac.in}




\begin{abstract}
Group/Pooled testing is a technique which involves testing $n$ samples for a rare disease in an indirect manner: instead of individually testing each of the $p$ samples, it involves creating $n < p$ pools and testing each pool, where a pool consists of a mixture of small, equal portions of a subset of the $p$ samples. Group testing is known to save a large amount of testing time and resources in a variety of application settings. There also exist good theoretical guarantees for the recovery of the status of the $p$ samples from results on $n$ pools, if certain conditions are satisfied. The popularity of group testing for RT-PCR testing is also well-known. The noise in quantitative RT-PCR is known to follow a multiplicative and data-dependent model, which emerges from properties inherent to RT-PCR. In recent literature, the corresponding linear systems for inference of the health status of $p$ samples from results on $n$ pools, have been solved using the well-known \textsc{Lasso} estimator and its variants \cite{Ghosh2021}, which have been widely used in the compressed sensing literature. However, the \textsc{Lasso} has typically been used in additive Gaussian noise settings. There is no existing literature which establishes performance bounds for \textsc{Lasso} for the multiplicative noise model associated with RT-PCR. After noting that a general technique proposed in \cite{Hunt2018} works well for establishing performance bounds for \textsc{Lasso} in Poisson inverse problems, we adapt it to handle sparse signal reconstruction from compressive measurements with multiplicative noise, as would appear in pooled RT-PCR testing. In particular, we present high probability performance bounds and data-dependent weights for the \textsc{Lasso} and a weighted version of the \textsc{Lasso} under multiplicative lognormal Noise. We also show numerical results on simulated pooled RT-PCR data to empirically validate our developed bounds.
\end{abstract}




\end{frontmatter}


\section{Introduction}
\label{sec:introduction}
The COVID-19 pandemic led to massive lockdowns worldwide, causing disruptions in industry, education and the well-being of people. Symptoms of this virus were often similar to those of the common cold, flu, etc. This, combined with the highly infectious nature of COVID-19, made early detection of the virus in individuals for subsequent isolation, one of the significant challenges faced by testing laboratories. 

RT-PCR (Reverse Transcription-Polymerase Chain Reaction) \cite{Jawerth2020} is a highly reliable method for detecting viruses, and was widely deployed during the COVID-19 pandemic. It involves obtaining a sample from the nose or throat of a person, chemically treating it to extract viral RNA (if present in the sample), and then converting the RNA to DNA through a process called reverse transcription. Short DNA fragments complementary to the viral DNA are added to the sample, which is then passed through a DNA amplification process. The presence of the virus is detected using fluorescent markers, and if the total fluorescence exceeds a certain threshold, the sample is considered positive for the virus. The number of cycles of amplification $C_t$ (cycle threshold) required to exceed this threshold is reported. A larger viral load will produce a larger amount of complementary DNA and thus produce the critical level of fluorescence earlier on in the RT-PCR process, leading to a lower $C_t$ value. Likewise, smaller viral loads are associated with higher $C_t$ values.

However, the low prevalence rate of COVID-19 means that most test samples turn out to be negative, leading to wastage of valuable resources, if every sample has to be tested. In an effort to reduce the cost and burden of RT-PCR testing, various algorithms have been proposed, that allow for the pooling of samples in order to achieve the same or slightly lower diagnostic accuracy with fewer tests \cite{PooledWiki2021}. One such algorithm is Tapestry, which was introduced in \cite{Ghosh2021}, and which combines ideas from combinatorial group testing and compressed sensing (CS). Other 
methods such as \cite{Shental2020,Heiderzadeh2021} also use interesting ideas from group testing (GT) and CS. Another approach, described in \cite{Goenka2021,Goenka2021contact}, uses additional \textit{side information} such as family membership or via contact tracing information. This additional information is shown to further improve the accuracy of the detection algorithms. 
All these algorithms serve two purposes: (\textit{i}) To reduce the number of testing kits required per person (which were scarce in many countries during COVID-19 surges), and (\textit{ii}) To increase testing throughput by parallelizing multiple pooled tests. 

Simple two-round group testing schemes such as Dorfman testing \cite{Dorfman1943} had already been applied in the field by several research labs \cite{israeli_news}, \cite{healthcare_europe} for COVID-19 testing. However Dorfman testing is a two-round procedure where the output of the first round acts as input to the second round. This increases testing time significantly, because the second round cannot begin before completion of the first round, and because one round of RT-PCR takes about 3-4 hours to complete. On the other hand, methods such as Tapestry \cite{Ghosh2021} are non-adaptive and work in a single round. This saves on testing time, which is very critical during the peak of a pandemic. In fact, some recent research \cite{Larremore2021} argues that testing time is even more important than accuracy at the peak of a pandemic.  
Tapestry is based on principles of CS. CS \cite{Donoho2006} is a mathematical framework that allows for the efficient recovery of sparse signals from a small number of carefully chosen linear measurements. Let $\boldsymbol{x^*} \in \mathbbm{R}^p$ be a signal with only $s \ll p$ non-zero elements at unknown indices. The signal $\boldsymbol{x^*}$ is observed in an indirect manner via linear measurements of the form $\boldsymbol{y} = \boldsymbol{Ax^*} + \boldsymbol{w}$ where $\boldsymbol{A} \in \mathbbm{R}^{n \times p}, n \ll p$ is a so called sensing matrix with fewer rows than columns,  $\boldsymbol{y} \in \mathbbm{R}^n$ is a measurement vector, and $\boldsymbol{w} \in \mathbbm{R}^n$ is an additive noise vector. Using simple convex optimization procedures, it has been shown that the sparse vector $\boldsymbol{x^*}$ can be recovered stably and robustly from $\boldsymbol{y}, \boldsymbol{A}$ provided two (related) conditions are met: (\textit{i}) The number of measurements should be at least $O(s \ln p)$, and (\textit{ii}) The matrix $\boldsymbol{A}$ should obey the property that its nullspace accommodates no signal with $s$ or fewer than $s$ non-zero elements, barring the zero vector. CS has found numerous applications in signal processing, image reconstruction, and other fields and has been particularly useful for sparse linear inverse problems, such as those encountered in pooled RT-PCR testing. Note that most traditional GT algorithms (see \cite{Aldridge2019} for a survey) treat both $\boldsymbol{y}$ and $\boldsymbol{x^*}$ as binary. However, CS produces a distinct advantage over these traditional GT algorithms: it allows us to determine \emph{quantitative} information such as viral load, over and above a mere binary indicator of whether or not the sample is diseased \cite{Shental2020,Ghosh2021,Goenka2021,Goenka2021contact}. 

The well known \textsc{Lasso} estimator in statistics \cite{THW2015} and its variants are popular methods for solving sparse linear inverse problems in CS. The estimator is given as follows:
\begin{equation}
    \boldsymbol{\hat{x}} := \textrm{argmin}_{\boldsymbol{x}} \|\boldsymbol{y}-\boldsymbol{Ax}\|^2 + \lambda \|\boldsymbol{x}\|_1,
\end{equation}
where $\lambda > 0$ is a regularization parameter. \textsc{Lasso} has been used in group testing, including the Tapestry algorithm \cite{Ghosh2021}. Its variants such as the group \textsc{Lasso} \cite{Huang2010} are useful methods to exploit \textit{side information} to enhance group testing performance \cite{Goenka2021,Goenka2021contact}. There is a rich literature on theoretical performance bounds for the \textsc{Lasso} under additive Gaussian \cite[Chapter 11]{THW2015} and Poisson noise \cite{Hunt2018,Li2018}. \cite{Hunt2018} also provides a novel analysis of a \emph{weighted} version of \textsc{Lasso}, deriving data-dependent weights and chalking out a generalized framework for handling data-dependent noise, with emphasis on the Poisson noise model which is widely used for optical imaging systems. However, the nature of RT-PCR testing introduces multiplicative noise (as we shall see in detail in Sec.~\ref{ssec:problemform}), whose effective variance scales directly with the values of the unknown vector $\bm{x^*}$. While the \textsc{Lasso} has been shown to numerically perform well in this setting \cite{Ghosh2021}, its theoretical bounds for this noise model have not been established in the literature to the best of our knowledge. It is important to have a better understanding of the theoretical performance guarantees of the \textsc{Lasso} under multiplicative noise to rigorously justify the use of this estimator in a multiplicative noise setting such as that encountered in pooled RT-PCR. 

In this paper, we draw from a general technique proposed in \cite{Negahban2012} and adapt the data-dependent performance bounds for \textsc{Lasso} in Poisson inverse problems as presented in \cite{Hunt2018}, to the case of multiplicative noise in the context of pooled RT-PCR testing. We derive high-probability performance bounds, obtained via data-dependent weights for \textsc{Lasso} and its weighted variant (details provided later in this paper), assuming randomly generated Bernoulli pooling matrices. We support our theory through numerical experiments. 

In parallel to this work, our group has done other work \cite{Rishav2021} on the analysis of the \textsc{Lasso} for multiplicative noise in compressive RT-PCR. This work may be submitted for publication later, independent of the work presented in this paper. There are two main differences between \cite{Rishav2021} and this work: (\textit{i}) the former does not explore a weighted version of \textsc{Lasso}, as we do here (see Equation \ref{eqn:wlasso}), and (\textit{ii}) unlike \cite{Rishav2021}, the work here considers a Taylor series based expansion of the multiplicative noise model (see Eqns.~\ref{eqn:cs_prob_noise} and \ref{eqn:noisemodel}).

Lastly, we note that non-Gaussian noise has been extensively considered in many branches of signal processing. Poisson noise is considered in compressed sensing in \cite{Rajwade2021,Bohra2019}. Moreover, \cite{Stojanovic2016} deals with optimal experiment design for ARX models under non-Gaussian noise. Likewise, \cite{Stojanovic2016_2} deals with robust identification of output error models in the presence of heavy-tailed, non-Gaussian noise. On the other hand, our paper considers a very specific non-Gaussian noise model -- multiplicative lognormal noise, which arises in RT-PCR. Moreover it deals with compressive recovery unlike the two aforementioned papers.

\section{Theory}
\label{sec:theory}

\subsection{Problem Formulation}
\label{ssec:problemform}
In RT-PCR testing, a naso- or oro-pharyngeal swab is taken from a person being tested. The mucus in the swab is mixed in a liquid medium and tested in an RT-PCR machine. \cite{surveygt} surveys multiple group testing techniques in detail, where they pool (mix) subsets of $p$ samples and test the $n$ pools ($\#\text{pools}(n) < \#\text{samples}(p)$) to save resources. Small, equal portions of participating samples are used for creating any pool. Given a negative pooled test, one can conclude that all contributing samples are non-infected, as RT-PCR is known to be highly sensitive by design due to the DNA amplification involved. In case some of the pooled tests are positive, we would like to deduce which participants are infected without repeated rounds of testing. This part is non-trivial. 

More formally, we can state the computational problem as follows. Let $\boldsymbol{x^*} \in \mathbbm{R}^{p}$ represent the vector of viral loads of the samples collected from $p$ individuals (one per individual). All entries of $\bm{x^*}$ are non-negative and $x_i = 0$ would imply that the sample of the $i$th individual is not infected. To model the low infection rate for COVID-19, we consider $\boldsymbol{x^*}$ to be a sparse vector with very few positive valued elements. The mixing of samples into groups (or pools) can be modeled with the help of a binary mixing matrix $\boldsymbol{A} \in \{0,1\}^{n \times p}$ (also called pooling matrix or sensing matrix), where $n < p$ and where entry $A_{ij}$ is 1 if the sample from the $j$th person is included while creating the $i$th pool and 0 otherwise. The resulting viral loads in pooled samples will be given by $\boldsymbol{y} \in \mathbbm{R}^{n}$, where:
$y_j = \sum_{i=1}^{p} A_{ji} x^*_i = \bm{A^{j}}\bm{x^*}, \, 1 \leq j \leq  n.$ Here, $\bm{A^{j}}$ is the $j$th row of $\bm{A}$.

However, due to the stochasticity of the reaction, measurement error in the PCR machine, and exponential growth of the number of molecules of viral cDNA in all the pools, a multiplicative lognormal noise model is used for $\boldsymbol{y}$. We use the noise model derived in \cite{Ghosh2021} for our analysis:
\begin{equation}\label{eqn:cs_prob_noise}
    y_j = \bm{A^j}\bm{x^*}(1 + q_a)^{w_j}.
\end{equation}
Here, $q_a \in (0,1)$ is the known fraction of viral cDNA that replicates in each cycle (sometimes referred to as the amplification factor), and $w_j \sim \mathcal{N}(0,\,\sigma^{2})$ is the sample of a random variable which represents the difference between the true cycle time and observed/recorded cycle time. Here $\sigma$ is the standard deviation of $w_j$. Usually, the time difference is small, and
hence we can assume $\sigma < 1$ (Fig.1 of Additional File 3 in \cite{Karlen2007}), as has also been done in \cite{Ghosh2021,Goenka2021}.

Several classes of binary pooling matrices have shown promising results when applied to the problem of inferring $\bm{x^*}$ from $\bm{y}, \bm{A}$, as shown in \cite{Ghosh2021}. This includes different types of deterministically generated matrices. However, these typically have restrictions on the number of rows and columns as discussed in \cite{Ghosh2021}. Hence, for greater flexibility, and to simplify our analysis, we limit our study to the class of random Bernoulli sensing matrices, where each element of the matrix is drawn independently from the $\textrm{Bernoulli}(q)$ distribution for some $q \in (0,1)$. That is, each entry $A_{ij}$ of matrix $\boldsymbol{A}$ is independently 0 with probability $1-q$ and 1 with probability $q$. Bernoulli matrices are useful for simplicity of pool creation: a sample either participates in a pool/group ($A_{ij} = 1$) or it does not ($A_{ij} = 0)$, whereas fractional entries would make pooling much more cumbersome. 

Our study is in the low variance regime, i.e. $\sigma < 1$, which we would later see appears to be an essential assumption for accurate or near-accurate recovery of $\boldsymbol{x^*}$ from $\boldsymbol{y}, \boldsymbol{A}$. This assumption is also supported in practice by RT-PCR systems, since the variance of the cycle time readings is very low as mentioned earlier. Additionally, this assumption allows us to approximate the multiplicative lognormal noise model by linearizing \cref{eqn:cs_prob_noise} using a first order Taylor series expansion, as also mentioned in \cite[Sec. VII]{Goenka2021}:
    \begin{equation}\label{eqn:noisemodel}
        y_j \approx \bm{A^j}\bm{x^*} + [\bm{A^j x^*}\ln(1+q_a)]w_j, \text{ where } w_j \sim \mathcal{N}(0, \sigma^2).
    \end{equation}
Our target is to recover vector $\bm{x^*}$ using our pooled viral loads $\bm{y}$ and the pooling matrix $\bm{A}$. Note that the approximation in the above equation is tight since $\sigma < 1$, all powers of $\ln (1+q_a)$ are less than 1 and due to division by factors such as $2!, 3!$, etc. All this makes the second-order or higher order terms of the Taylor series expansion negligible in magnitude. The Taylor expansion of $y_j$ is as follows:
    \begin{equation}
        y_j \approx \bm{A^j}\bm{x^*}\left( 1 + \ln(1+q_a)w_j + \frac{1}{2}\ln^2(1+q_a)w_j^2 + \cdots \right)
    \end{equation}
As $w_j \sim \mathcal{N}(0,\sigma^2)$ and $\sigma < 1$, we can ignore the higher order terms of $w_j^2$, $w_j^3$, $\cdots$, and end up with \cref{eqn:noisemodel}. Ideally, this approximation is valid only with high probability. But for simplicity, we treat it as exact. 

The Taylor series expansion casts this model as a special case of additive signal dependent noise. Note that the additive noise variance in this model (see \cref{eqn:noisemodel}) is proportional to the \emph{square} of the mean, unlike the case of Poisson noise, where the noise variance is \emph{equal} to the mean. Along a related vein, the work in \cite{Arildsen2014} presents an interesting case where the noise in the compressive measurements is linearly correlated with the signal, and many applications including quantization noise are handled. However no theoretical performance bounds are presented, and multiplicative noise is not considered. 


\subsection{Estimation methods}\label{ssec:estmeth}
We consider two estimators in this work: the \textsc{Lasso} and the weighted \textsc{Lasso} (referred to henceforth as \textsc{WLasso}), both for the case of multiplicative noise in $\bm{y}$. 

The \textsc{Lasso} estimator is presented here below:
\begin{equation}\label{eqn:lasso}
    \bm{\hat{x}^L} = \operatorname*{argmin}_{\bm{x} \in \mathbbm{R}^{p}} \|\bm{y} - \bm{Ax}\|_{2}^{2} + \gamma \sum_{k=1}^{p}\beta|x_{k}|,
\end{equation}
where $\gamma > 2$ is a constant and $\beta > 0$ is a data-dependent scalar.

The \textsc{WLasso} estimator uses a possibly different weight for each element of $\bm{x}$, as shown below:
\begin{equation}\label{eqn:wlasso}
    \bm{\hat{x}^{WL}} = \operatorname*{argmin}_{\bm{x} \in \mathbbm{R}^{p}} \|\bm{y} - \bm{Ax}\|_{2}^{2} + \gamma \sum_{k=1}^{p}\beta_k |x_{k}|,
\end{equation}
where $\gamma > 2$ is a constant and $\beta_k > 0$ is the $k$th data-dependent scalar. The exact expressions of $\beta$ and $\beta_k$ depend on the observed data. Deriving expressions for $\beta$ and $\beta_k$ forms the crux of our work in this paper (albeit for modified forms of these estimators presented later in \cref{eqn:lassosurr} and \cref{eqn:wlassosurr}).

Additionally, note that a logarithmic transformation on $\bm{y}$ will act as a variance stabilization transform \cite{Pollard1979} for multiplicative noise, just as the square-root transform acts a variance stabilizer for Poisson noise and has shown great success in Poisson/Poisson-Gaussian inverse problems \cite{Rajwade2021,Bohra2019}, tomography \cite{Gopal2020} and deblurring \cite{Dupe2009}. The associated estimator for our multiplicative noise outlined in \cref{eqn:cs_prob_noise}, i.e. 
 $\textrm{argmin}_{\bm{x}}\lVert \ln \bm{y} - \ln \bm{Ax} \rVert^2_2 + \lambda \| \bm{x} \|_1$, might seem like an obvious choice. However, we have observed that such an estimator produces very poor results in practice, and has many convergence issues. Hence, we have not pursued it further in this work. Similar observations have also been reported in \cite[Sec. VII-A]{Goenka2021}. 

 Let $S^{*}$ be the support of $\bm{x^{*}}$ (the true, unknown vector of $p$ viral loads) having size $s = |S^{*}|$. Under this condition of exact sparsity, \cite[Corollary 2]{Negahban2012} states that if the sensing matrix $\bm{A}$ satisfies the so-called Restricted Eigenvalue Condition (REC, defined later) and the column normalization condition, then the optimal solution $\bm{\hat{x}^L}$ of the \textsc{Lasso} instance, with regularization parameter $\beta=4\sigma\sqrt{\frac{\ln p}{n}}$, satisfies a high probability performance bound which increases in direct proportion to the regularization parameter. However, \cite{Negahban2012} sets a requirement on $\beta$ for this bound to hold: $\beta\geq 2\|\nabla \mathcal{L}\|_{\infty}$ where $\mathcal{L}(.)$ is the likelihood function, i.e., $\beta$ should be greater than twice the $l_{\infty}$ norm of the gradient of the quadratic loss in \textsc{Lasso}. This becomes particularly restrictive when we apply it to the general case of \textsc{WLasso} from \cref{eqn:wlasso}, where we have a separate regularization parameter $\beta_k$ for each element of the signal. Instead, we can extend this by lower bounding each $\beta_k$ by the absolute value of its individual gradient term, i.e. by $|(\bm{A}^T(\bm{y} - \bm{Ax^*}))_k|$. 

Under this new condition on $\beta_k$, we have the following risk bounds for the weighted and classical \textsc{Lasso} estimates ($\bm{\hat{x}^{WL}}$ and $\bm{\hat{x}^L}$ respectively) when $\bm{A}$ and $\bm{y}$ follow the conditions of \cref{prop3} outlined later in the paper. Here $\rho_{\gamma}$ only depends on $\gamma$ and $\eta$ is a parameter associated with the restricted eigenvalue condition (REC) of the sensing matrix $\bm{A}$:
\begin{equation}
\|\bm{\hat{x}^{WL}}-\bm{x^{*}}\|_{2}^{2} \leq \frac{\rho_{\gamma}^2}{\eta}\sum_{k \in S^{*}} \beta_k^2; \textrm{  } \|\bm{\hat{x}^{L}}-\bm{x^{*}}\|_{2}^{2} \leq \frac{\rho_{\gamma}^2}{\eta} s \beta^2.
\label{eq:lasso_bounds}
\end{equation}
Looser bounds than those shown here can be obtained from an analysis of the \textsc{WLasso} estimator using the typical \textsc{Lasso} bounding procedures outlined in \cite{Negahban2012}. Following from the prior discussion, if each $\beta_k$ in the \textsc{WLasso} estimator is close to $|(\bm{A}^T(\bm{y} - \bm{Ax^*}))_{k}|$, then the bounds on the \textsc{WLasso} are close to those of the oracle least squares estimator \cite{Negahban2012,Hunt2018}. But practically, we cannot set $\beta_k$ to these values as the weights should only depend on the observed data. However, we can still obtain reasonable estimates of these weights based on the characteristics of the data, as we do in this work. 
\begin{equation}\label{eqn:asswts}
\textrm{Condition } \mathscr{C}1 \textrm{ on } \left(\{\beta_k\}_{k}\right):  |(\bm{A}^T(\bm{y} - \bm{Ax^*}))_{k}| \leq \beta_k \quad \forall  k=1,...,p.
\end{equation}
We want the \emph{smallest} possible values of $\beta_k$ which would faithfully satisfy the Condition $\mathscr{C}1$ with high probability. This serves as our primary objective in deriving the data-dependent weights for our noise model as stated in \cref{ssec:assandddw}. 

\subsection{Assumptions and Data-Dependent Weights}\label{ssec:assandddw}
Here we present the data-dependent weights and the underlying assumptions we propose for the \textsc{Lasso} and \textsc{WLasso} estimators defined in Equation \ref{eqn:lasso} and Equation \ref{eqn:wlasso}, respectively. In this section, we outline the overall approach and a broad proof sketch, while a full proof for these expressions can be found in the Supplemental Material. We further demonstrate simulation results in \cref{sec:numerical}, and theoretical performance bounds in \cref{ssec:perfbound} using these weights. 

We first impose a condition on the sensing matrix required for reconstruction using the variants of the \textsc{Lasso} estimator. This is known as the \textit{Restricted Eigenvalue Condition}. It is defined below here, following \cite[Chapter 11]{THW2015}: 

\noindent A sensing matrix $\boldsymbol{A}$ of size $n \times p$ is said to obey the Restricted Eigenvalue Condition (REC), if there exist constants $\kappa_1, \kappa_2 > 0$ such that
\begin{equation}
    \lVert \bm{Ax} \rVert_2 \geq \kappa_2 \lVert \bm{x} \rVert_2 - \kappa_1 \lVert \bm{x} \rVert_1, \quad \forall \bm{x} \in \mathbbm{R}^{p}.
\label{eqn:rec} 
\end{equation}

\subsubsection{Rescaling and Recentering}\label{ssec:rescalerecenter}
The random Bernoulli sensing matrix which is suitable for pooled testing does not obey the REC. Instead, we define a surrogate pooling matrix $\bm{\Tilde{A}}$ which would obey the REC and a surrogate measurement vector $\bm{\Tilde{y}}$. We want a complementary transformation to generate $\bm{\Tilde{y}}, \bm{\Tilde{A}}$ from $\bm{y}, \bm{A}$. Ideally, we want the Gram matrix of our surrogate pooling matrix to be as homogeneous as possible. We invoke the following two commonly used conditions to obtain our surrogate matrices:
\begin{equation}\label{eqn:surrogateA}
  \mathbbm{E}\left(\bm{\Tilde{A}}^{T}\bm{\Tilde{A}}\right) = \bm{I_p}; 
    \mathbbm{E}\left(\bm{\Tilde{A}}^T\left(\bm{\Tilde{y}} - \bm{\Tilde{A}x^*}\right)\right) = \bm{0}.
  \end{equation}
In \cref{eqn:surrogateA}, $\bm{I_p}$ stands for a $p \times p$ identity matrix. The first constraint in \cref{eqn:surrogateA} helps us ensure that the REC defined in \cref{eqn:rec} holds. We impose the second constraint in \cref{eqn:surrogateA} to make $\beta_k$ as small as possible, while satisfying the Condition $\mathscr{C}1$. Based on the calculations in \cite[Appendix C-A]{Hunt2018}, we obtain the following surrogate sensing matrix $\bm{\Tilde{A}}$ and surrogate measurement vector $\bm{\Tilde{y}}$ given a Bernoulli$(q)$ sensing matrix $\bm{A}$ and its associated measurement vector $\bm{y}$:
\begin{equation}
    \bm{\Tilde{A}} = \frac{\bm{A}}{\sqrt{nq(1-q)}} - \frac{q\bm{1}_{n\times 1}\bm{1}_{p\times 1}^T}{\sqrt{nq(1-q)}}; 
    \bm{\Tilde{y}} = \frac{1}{(n-1)\sqrt{nq(1-q)}}\left( n\bm{y} - \sum_{l=1}^{n}y_l\bm{1}_{n\times 1} \right).
\end{equation}
We use the surrogate matrix $\boldsymbol{\tilde{A}}$ and surrogate measurement vector $\boldsymbol{\tilde{y}}$ within the \textsc{Lasso} and \textsc{WLasso} estimators, as will be shown in the next section. We note that $\boldsymbol{\tilde{A}}$ and $\boldsymbol{\tilde{y}}$ can be computed entirely from $\boldsymbol{y}, \boldsymbol{A}, q$ without knowledge of the underlying signal. Moreover, we note that the REC holds with high probability for $\boldsymbol{\tilde{A}}$, as established in the lemma below which was proved in \cite[Appendix C-B]{Hunt2018}.
\begin{lemma} 
\label{lemma1}
There exist positive constants $c, c', c''$ such that with probability larger than $1-c'\exp(-c''n)$,
$\boldsymbol{\tilde{A}}$ satisfies REC($\kappa_1, \kappa_2$) with
parameters $\kappa_1 := \dfrac{c}{q(1-q)}\sqrt{\dfrac{\ln p}{n}}, \kappa_2 := \dfrac{1}{4}$.
\end{lemma}
\noindent

\subsubsection{Data dependent weights}
Given the surrogate sensing matrix $\bm{\Tilde{A}}$ and observations $\bm{\Tilde{y}}$, we redefine the \textsc{Lasso} and \textsc{WLasso} estimators via Equations \eqref{eqn:lassosurr} and \eqref{eqn:wlassosurr} respectively, as given below: 
\begin{equation}\label{eqn:lassosurr}
    \bm{\hat{x}^L} = \operatorname*{argmin}_{\bm{x} \in \mathbbm{R}^{p}} \|\bm{\Tilde{y}} - \bm{\Tilde{A}x}\|_{2}^{2} + \gamma \sum_{k=1}^{p}\beta|x_{k}|,
\end{equation}
where $\gamma > 2$ is a constant and $\beta > 0$ is a data-dependent scalar.
\begin{equation}\label{eqn:wlassosurr}
    \bm{\hat{x}^{WL}} = \operatorname*{argmin}_{\bm{x} \in \mathbbm{R}^{p}} \|\bm{\Tilde{y}} - \bm{\Tilde{A}x}\|_{2}^{2} + \gamma \sum_{k=1}^{p}\beta_k |x_{k}|,
\end{equation}
where $\gamma > 2$ is a constant and $\beta_k > 0$ is the $k$th data-dependent scalar. The exact expressions of $\beta$ and $\beta_k$ depend on the observed data and form the crux of our work.

In this section, we now upper bound the term $\left|\left(\bm{\Tilde{A}}^T\left(\bm{\Tilde{y}} - \bm{\Tilde{A}x^*}\right)\right)_{k}\right|$ and present data-dependent weights which satisfy Condition $\mathscr{C}1$ with high-probability. We derive a single weight $\beta$ for \textsc{Lasso} in \cref{prop1}, and then derive the weights $\{\beta_k\}_{k=1}^p$ for the \textsc{WLasso} in \cref{prop2}. The proofs of both propositions can be found in the Supplemental Material. 

\medskip
\begin{proposition} 
{\normalfont\textbf{Weight for \textsc{Lasso}:}}
\label{prop1}

\noindent
Consider $\bm{A},\bm{\bar{R}} \in \mathbbm{R}^{n \times p}$ and define $W := \max_{i,j} \left(\bm{\bar{R}}^T \bm{A} \right)_{ij}$ for $1 \leq i \leq p, 1 \leq j \leq p, \kappa := \sigma \ln (1+q_a)$, where $\bm{\bar{R}}_{k,l} := \left(\frac{na_{l,k} - \sum_{l'=1}^{n}a_{l',k}}{n(n-1)q(1-q)}\right)^2$ for $k = 1, \cdots, p; l = 1,\cdots, n$. Consider the following assumptions:
\begin{center}
\begin{itemize}
\addtolength{\itemindent}{1em}
    \item[\textit{A1.}] $nq \geq 12\max(q,1-q)\ln(p)$
    \item[\textit{A2.}] $p \geq 2$
    \item[\textit{A3.}] $\sigma < \frac{1}{\sqrt{2}\ln \left( 1 + q_a \right)} \left(\ln\left(\frac{1}{1 - \left( 1 - \frac{1}{p^3} \right)^{1/n}} \right)\right)^{-1/2}$
    \item[\textit{A4.}] $\bm{A}$ \textit{is a Bernoulli Sensing Matrix}.
\end{itemize}    
\end{center}
\noindent
If the aforementioned assumptions hold, then there exist positive constants $c$, $c'$ such that with probability larger than $1 - \frac{c'}{p}$, the choice%
\begin{equation*}
    \beta := \kappa \hat{\Lambda} \sqrt{6W\ln(p)} + c\left( \frac{3\ln(p)}{n} + \frac{9\max(q^2, (1-q)^2)}{n^2q(1-q)}\ln(p)^2 \right) \hat{\Lambda}
\end{equation*}%
satisfies Condition $\mathscr{C}1$ (defined earlier in \cref{eqn:asswts}), where $\hat{\Lambda}$ is an estimator of $\lVert \bm{x^*} \rVert_1$ given by
\begin{equation*}
    \hat{\Lambda} := \frac{\sum_{i=1}^{n} y_i + \sqrt{\sum_{i=1}^{n} y_i^2} \frac{\kappa\sqrt{6\ln(p)}}{1 - \kappa\sqrt{2g(3\ln(p))}}}{nq - \sqrt{6nq(1-q)\ln(p)} - \max(q, 1-q)\ln(p)},
\end{equation*}
where for any $\theta \in \mathbb{R}$, we define
$    g(\theta) := \ln\left( \frac{1}{1 - \left( 1 - e^{-\theta} \right)^{1/n}} \right)$.
Furthermore, $c=126$ works as long as $n \geq 20$.
\end{proposition}

\bigskip
\begin{proposition} {\normalfont\textbf{ Weights for \textsc{WLasso}}}\label{prop2}

\noindent
With the same notations and assumptions as in Proposition 1, there exist positive constants $c$, $c'$ such that with probability larger than $1 - \frac{c'}{p}$, the choice (depending on $k$)
\begin{equation*}
    \beta_k := \sqrt{\bm{\bar{R}_{k}}^T \bm{y_2}}\frac{\kappa\sqrt{6\ln(p)}}{1 - \kappa\sqrt{2g(3\ln(p))}} + c\left( \frac{3\ln(p)}{n} + \frac{9\max(q^2, (1-q)^2)}{n^2q(1-q)}\ln(p)^2 \right) \hat{\Lambda}
\end{equation*}
\noindent
satisfies Condition $\mathscr{C}1$ (see \cref{eqn:asswts}), where $\bm{y_2}:=\bm{y}\odot\bm{y}$. Furthermore, $c=126$ works as long as $n \geq 20$.
\end{proposition}

\bigskip
\noindent
\textbf{Comments on \cref{prop1} and \cref{prop2}: }

\begin{enumerate}

    \item We now argue that Assumption \textit{A3} is reasonable for our RT-PCR noise model. Since $n \geq 20$ and $p > n$, we have $1-\left(1-\frac{1}{p^3}\right)^{\frac{1}{n}} > \frac{1}{np^3}$ by Bernoulli's inequality \cite{BernoulliInequality_wiki}, \cite[page 20]{Brannan2006}. Due to this, the expression $\frac{1}{\sqrt{2}\ln \left( 1 + q_a \right)}\ln^{-\frac{1}{2}}\left( np^3 \right)$ acts as a lower bound on the RHS of the inequality in Assumption \textit{A3} of \cref{prop1}. This lower bound decreases as $n$ increases, as well as when $q_a$ increases. Since we are operating in a compressive regime, $n$ can be no greater than $p$. Moreover the largest possible value of $q_a$ is 1. Hence the lower bound on the LHS is given by $\dfrac{1}{2\ln(2)\sqrt{2\ln(p)}}$. Therefore, if we consider
$\sigma < \dfrac{1}{2\ln(2)\sqrt{2\ln(p)}}$, we always satisfy Assumption \textit{A3} of \cref{prop1}. Given the small values of $\sigma$ required in practice for the RT-PCR setting, as well as for the first-order expansion of the multiplicative noise model described in \cref{eqn:noisemodel}, this restriction on $\sigma$ is reasonable and holds true for a large range of $p$. For example, if $p=100$, we have $\sigma \leq 0.2377$.
      
\item Under Assumption \textit{A1} of \cref{prop1}, we avoid values of $q$ very close to 0 or 1. Hence the weights $\beta$ or $\{\beta_k\}_{k=1}^p$ are not very large in value, as can be seen from \cref{prop1} and \cref{prop2}. Thereby the upper bounds of the \textsc{Lasso} and \textsc{WLasso} estimators are also not very large, as will be established in Proposition \ref{prop3}.
    

\item With Assumptions \textit{A1} and \textit{A3} in place, we can show that $\hat{\Lambda}$, that is the estimator for $\|\bm{x}^*\|_1$, is always positive. Note that Assumption \textit{A3} can be re-written as $\kappa\sqrt{2g(3\ln p)} < 1$ by setting $\theta = 3 \ln p$ and using the definition of $\kappa$. This makes the numerator of the expression for $\hat{\Lambda}$ positive, whenever \textit{A3} is true as each element of $y$ is non-negative. We now argue that the denominator is also positive, for which we denote $m_q := \textrm{max}(q,1-q)$ and set $nq = h m_q \ln p$ for $h \geq 12$ from Assumption \textit{A1}. Then the denominator becomes $\left((h-1)m_q  - \sqrt{6(1-q)hm_q}\right)\ln p$ which is greater than or equal to $\left((h-1)m_q  - m_q\sqrt{6h}\right)\ln p$ using the definition of $m_q$. For any $h \geq 8$, we see that this expression is always positive, since $m_q$ and $\ln p$ are both positive. At this point, we have established only lower bounds on $\hat{\Lambda}$ and have not proved any upper bounds on its value in terms of $\|\bm{x^*}\|_1$. However, as will be shown in Sec.~\ref{sec:numerical}, the estimates are close to $\|\bm{x^*}\|_1$ within a factor of two.
\end{enumerate}

\subsection{Performance Bounds}
\label{ssec:perfbound}
In this section, we establish performance bounds for \textsc{Lasso} and \textsc{WLasso} in the context of multiplicative lognormal Noise. For this, we present \cref{prop3}, which is derived directly from \cite[Proposition 1]{Hunt2018} that accounts for the weighted $l_1$ regularizer defined in \cref{eqn:wlassosurr}, combined with Propositions \cref{prop1} and \cref{prop2} derived here for our specific noise model. 

\begin{proposition}\label{prop3} {\normalfont\textbf{Performance Bounds for \textsc{Lasso} and \textsc{WLasso}}}

\noindent
Fix $\varepsilon > 0$ and consider the following:
\begin{gather*}
    \beta_{max} := \max_{k \in \{1,\cdots,p\}} \beta_k,     \beta_{min} := \min_{k \in \{1,\cdots,p\}} \beta_k,
    \rho_{\gamma} := \gamma\frac{\gamma+2}{\gamma-2} \\
    \bm{D} := diag(\beta_1, \cdots, \beta_k, \cdots, \beta_p),
\end{gather*}%
\noindent
where $\gamma, \beta, \{\beta_k\}$ are parameters of the estimators presented in \cref{ssec:estmeth}, and $\bm{D}$ is a diagonal matrix with the $k$th diagonal element equal to $\beta_k$. For any vector $\bm{z} \in \mathbbm{R}^p$ and any set $S \subseteq \{1, \cdots, p\}$, define vector $\bm{z_{S}}$ such $z_S(i) = z_i$ for all $i \in S$ and $z_S(i) = 0$ for all $i \notin S$.

\noindent
Now consider that the following assumptions hold true:
\begin{itemize}
\addtolength{\itemindent}{1em}
    \item[B1.] $\gamma > 2$
    \item[B2.] $\{\beta_k\}_{k=1}^p$ and $\beta$ satisfy the Condition $\mathscr{C}1$ stated in \cref{ssec:estmeth}
    \item[B3.] $\bm{\Tilde{A}}$ satisfies the \text{REC} $(\kappa_1, \kappa_2)$ property stated in \cref{eqn:rec}.
\end{itemize}
\noindent
Then, there exists a universal constant $c>0$ such that for any set $S \subseteq \{1, \cdots, p\}$ for which
\begin{equation*}
    \lVert \bm{\beta_S} \rVert_2 \leq \beta_{min}\frac{\kappa_2 - \varepsilon}{\kappa_1 \rho_{\gamma}},
\end{equation*}
the \textsc{WLasso} estimator satisfies
\begin{equation*}
    \lVert \bm{x^*} - \bm{\hat{x}^{WL}} \rVert_2 \leq \frac{c\rho_{\gamma}}{\varepsilon^2}\lVert \bm{\beta_S}\rVert_2. 
\end{equation*}
\noindent
Furthermore, for any set $S \subseteq \{1,\ldots,p\}$ with $s = |S|$ satisfying
\begin{equation}
    \sqrt{s} \leq \frac{\kappa_2 - \varepsilon}{\kappa_1 \rho_{\gamma}},
\end{equation}
\noindent
for $\varepsilon>0$, the \textsc{Lasso} estimator satisfies
\begin{equation}
    \lVert \bm{x^{*}} - \bm{\hat{x}^L} \rVert_2 \leq \frac{c\rho_{\gamma}}{\varepsilon^2} \beta\sqrt{s}. 
\end{equation}
\end{proposition}

\noindent
\textbf{Comments on \cref{prop3}:}
\begin{enumerate}
\item Performance bounds for \textsc{Lasso} in \cite[Exercise 11.3, Eqn. 11.44a]{THW2015} or \cite[Proposition 1]{Hunt2018} are derived for weakly sparse signals. In contrast, we operate in the purely sparse regime for $\bm{x^*}$, as weakly sparse signal models are not required for viral loads. 

\item Consider the performance bounds for $\bm{\hat{x}^{WL}}$ and $\bm{\hat{x}^L}$ and the Assumption \textit{B2} from \cref{prop3}. Taking these into account, we observe that our motivation in \cref{prop1} and \cref{prop2} to find the smallest possible weights which satisfies the Condition $\mathscr{C}1$ is well justified. Ideally we would want weights $\beta_k$ such that $\beta_k = |(\bm{A}^T(\bm{y} - \bm{Ax^*}))_{k}| \quad \forall  k=1,...,p$. However, since the input signal $\bm{x^*}$ is not observable, we can at best produce high probability approximations, as done in \cref{prop1} and \cref{prop2}.
    
\item Upper bounds for both estimators increase with the $s$, the $\ell_0$ norm of $\boldsymbol{x^*}$. As $n$ increases, the upper bounds decrease, because the terms $W$ in \cref{prop1}, and $\sqrt{\boldsymbol{\bar{R}_k}^T \boldsymbol{y_2}}$ as well as $\dfrac{3c\ln p}{n}$ in \cref{prop2}, all decrease with $n$. These trends with respect to $s$ and $n$ are quite intuitive. 

\item As $\|\boldsymbol{x^*}\|_1$ increases, we observe from \cref{prop1} that all the $y_i$ values will typically increase and hence $\hat{\Lambda}$ will increase. Due to this, the upper bounds on $\|\boldsymbol{x^*}-\boldsymbol{\hat{x}^{L}}\|_2$ and $\|\boldsymbol{x^*}-\boldsymbol{\hat{x}^{WL}}\|_2$ for \textsc{Lasso} and \textsc{WLasso} respectively will also increase. However, the \emph{relative} errors $\dfrac{\|\boldsymbol{x^*}-\boldsymbol{\hat{x}^{WL}}\|_2}{\|\boldsymbol{x^*}\|_1}$
and $\dfrac{\|\boldsymbol{x^*}-\boldsymbol{\hat{x}^L}\|_2}{\|\boldsymbol{x^*}\|_1}$ remain constant. We emphasize the difference between this phenomenon and the trends observed in Poisson compressed sensing \cite{Rajwade2021,Hunt2018,Li2018,Bohra2019}, where the relative error actually \emph{decreases} with $\|\boldsymbol{x^*}\|_1$. The reason for this difference is that the noise standard deviation is proportional to the \emph{square root} of the signal value in case of Poisson noise, and proportional to the signal value in the multiplicative noise model of RT-PCR (see \cref{eqn:noisemodel}). 
\end{enumerate}

\section{Numerical Results}
\label{sec:numerical}

\begin{figure}[t!]
    \centering 
\begin{subfigure}{0.33\textwidth}
  \includegraphics[width=1.1\linewidth]{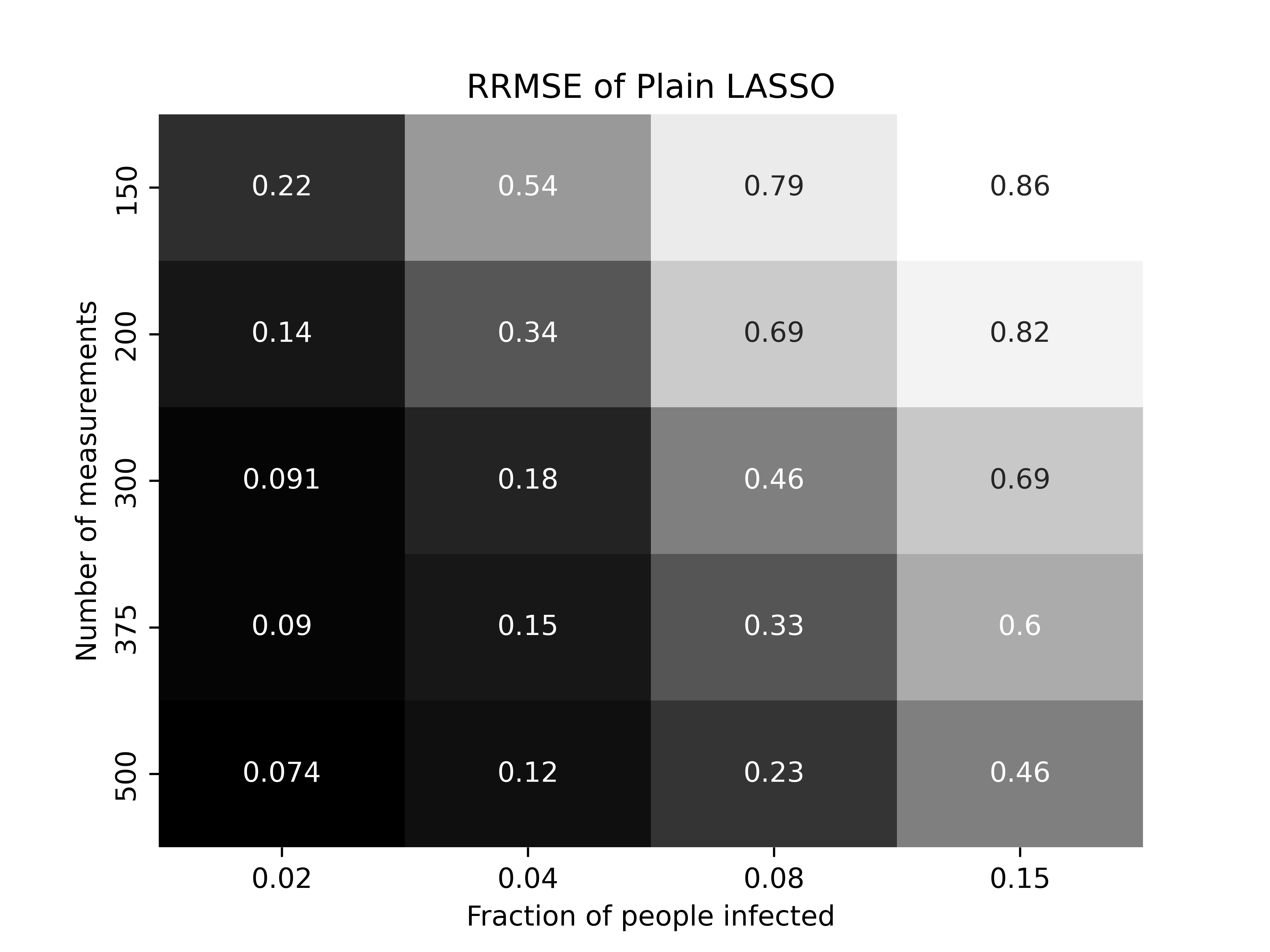}
  \label{fig:lasso_rrmse_0.5}
\end{subfigure}\hfil 
\begin{subfigure}{0.33\textwidth}
  \includegraphics[width=1.1\linewidth]{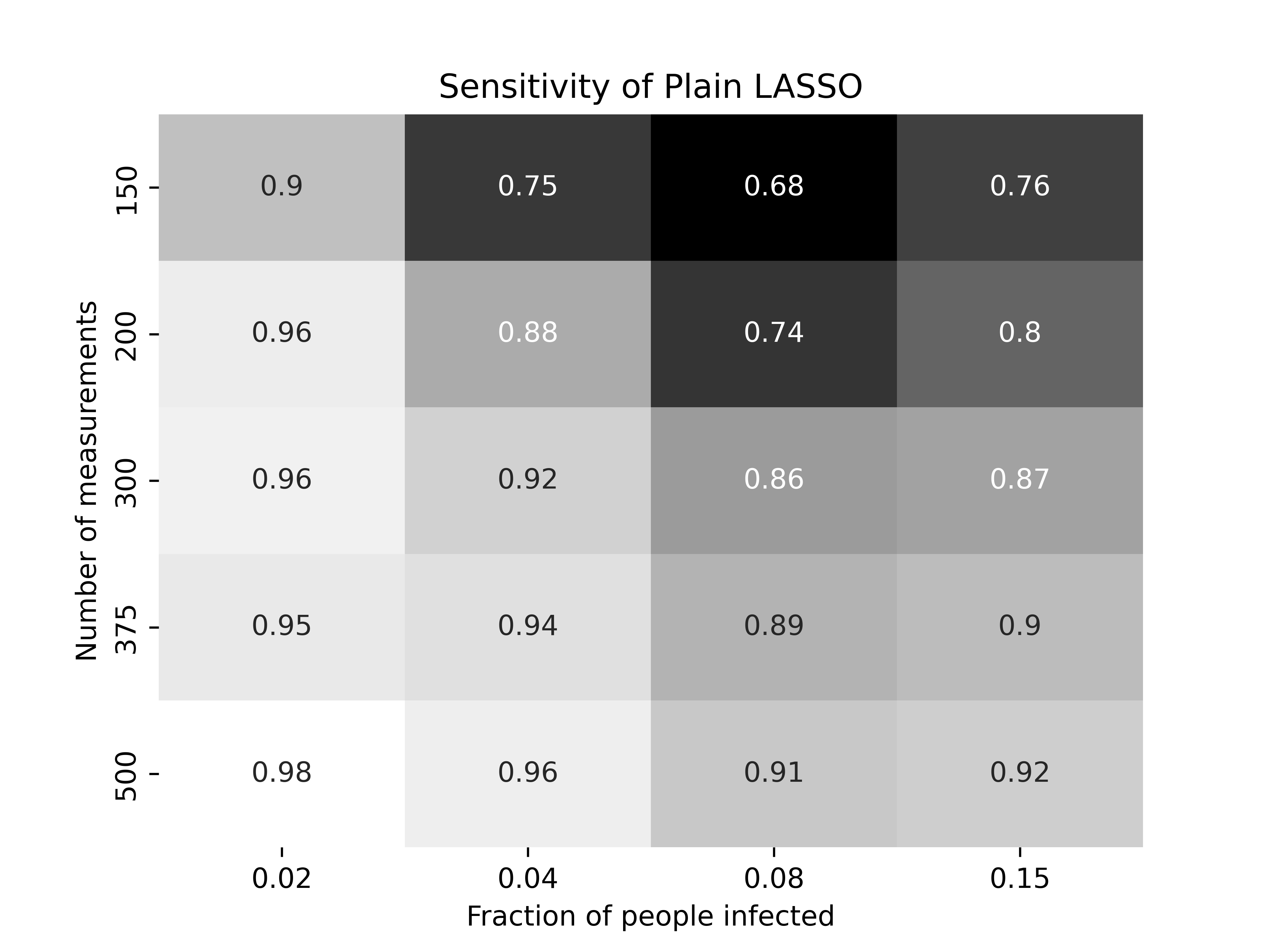}
  \label{fig:lasso_sen_0.5}
\end{subfigure}\hfil 
\begin{subfigure}{0.33\textwidth}
  \includegraphics[width=1.1\linewidth]{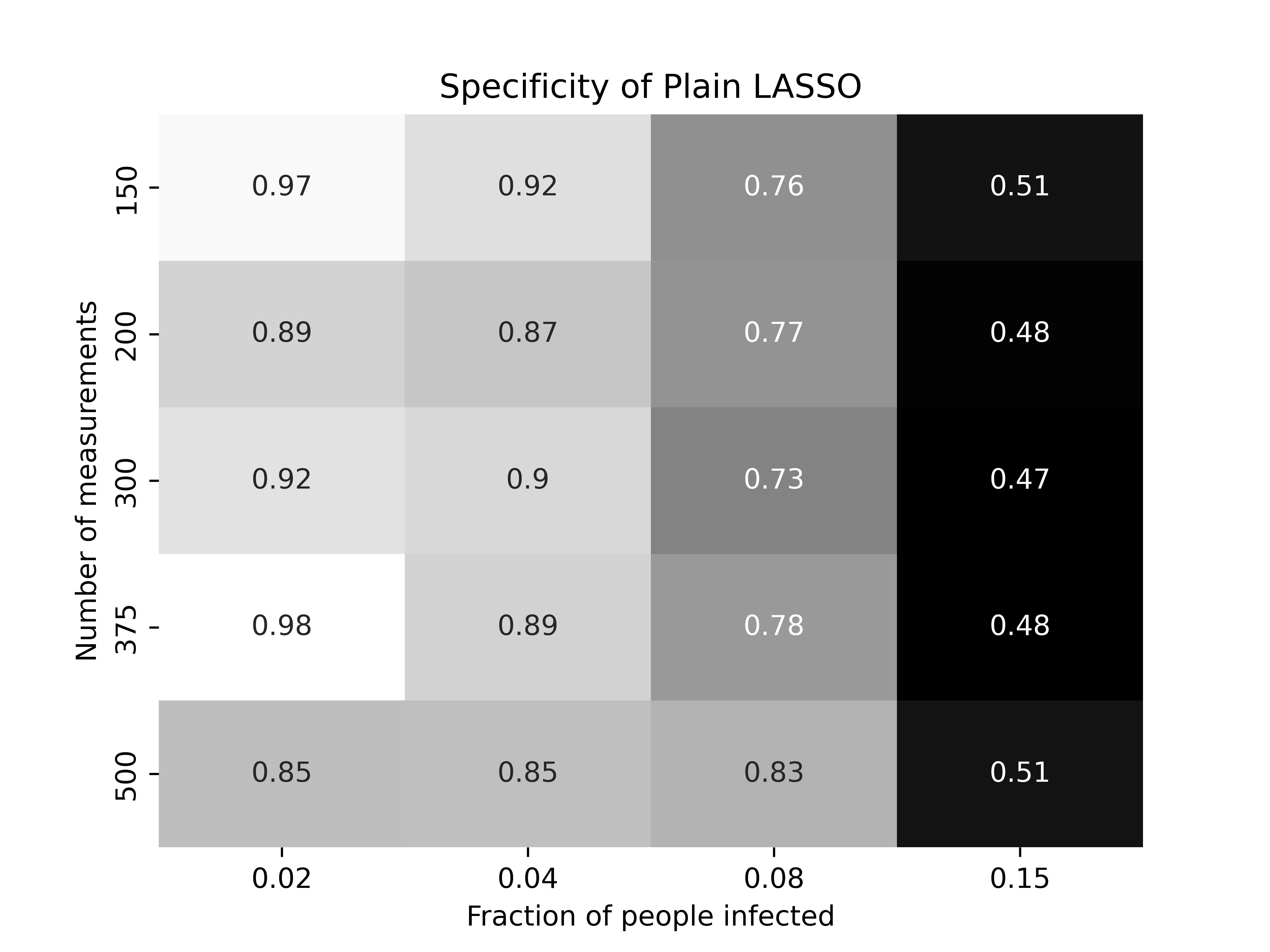}
  \label{fig:lasso_spec_0.5}
\end{subfigure}

\vspace{-1\baselineskip}
\begin{subfigure}{0.33\textwidth}
  \includegraphics[width=1.1\linewidth]{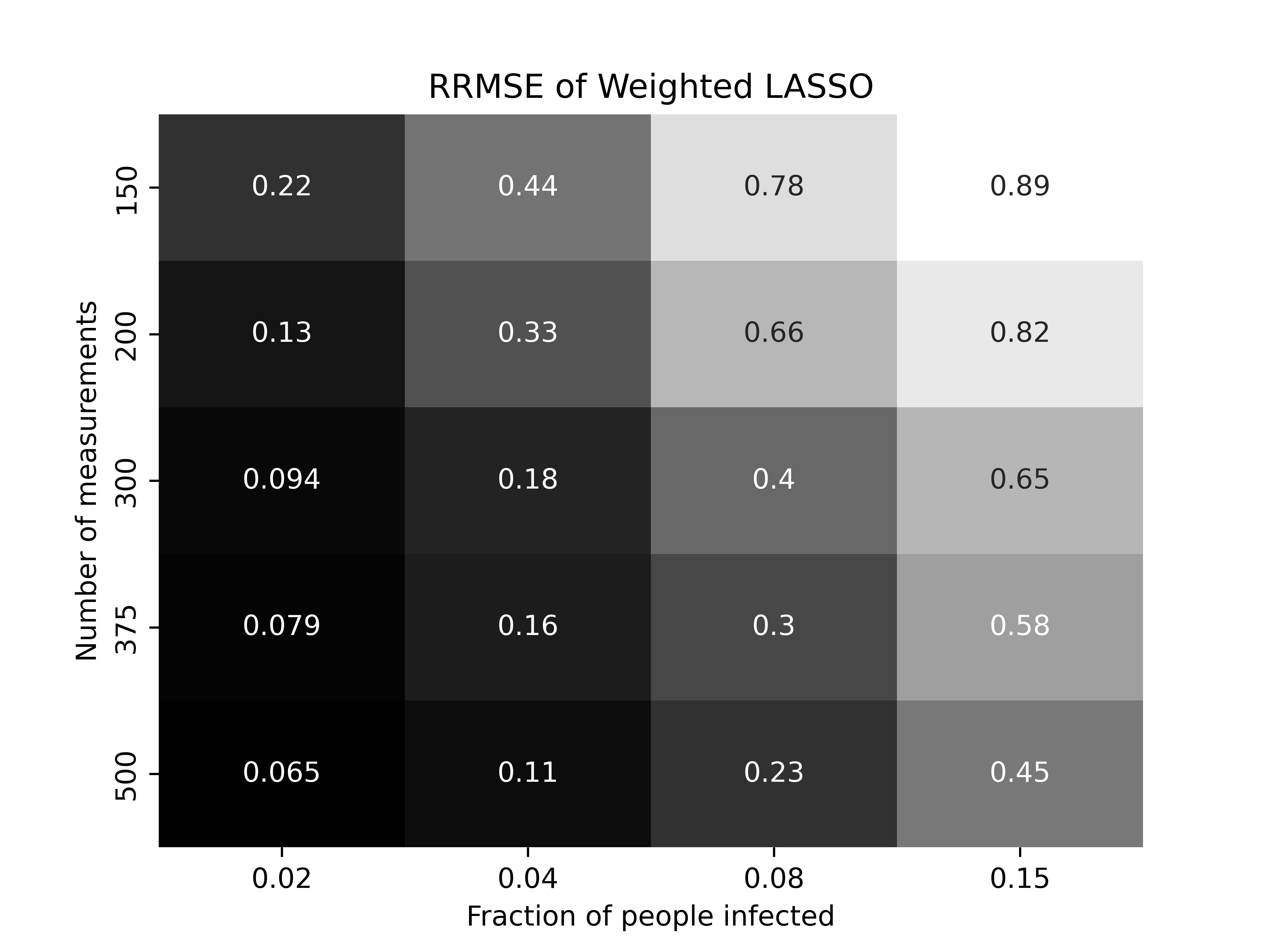}
  \label{fig:wlasso_rrmse_0.5}
\end{subfigure}\hfil
\begin{subfigure}{0.33\textwidth}
  \includegraphics[width=1.1\linewidth]{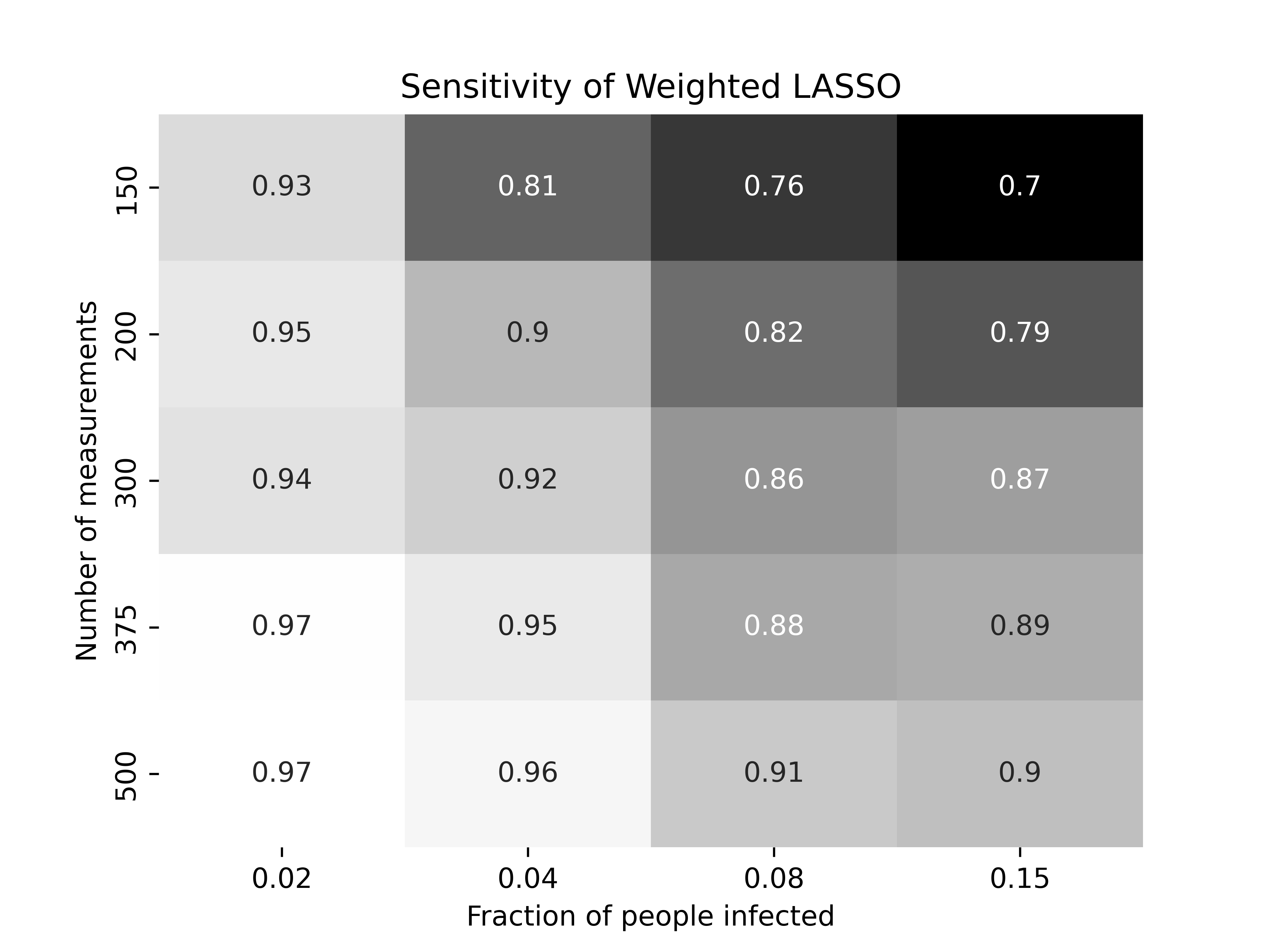}
  \label{fig:wlasso_sen_0.5}
\end{subfigure}\hfil
\begin{subfigure}{0.33\textwidth}
  \includegraphics[width=1.1\linewidth]{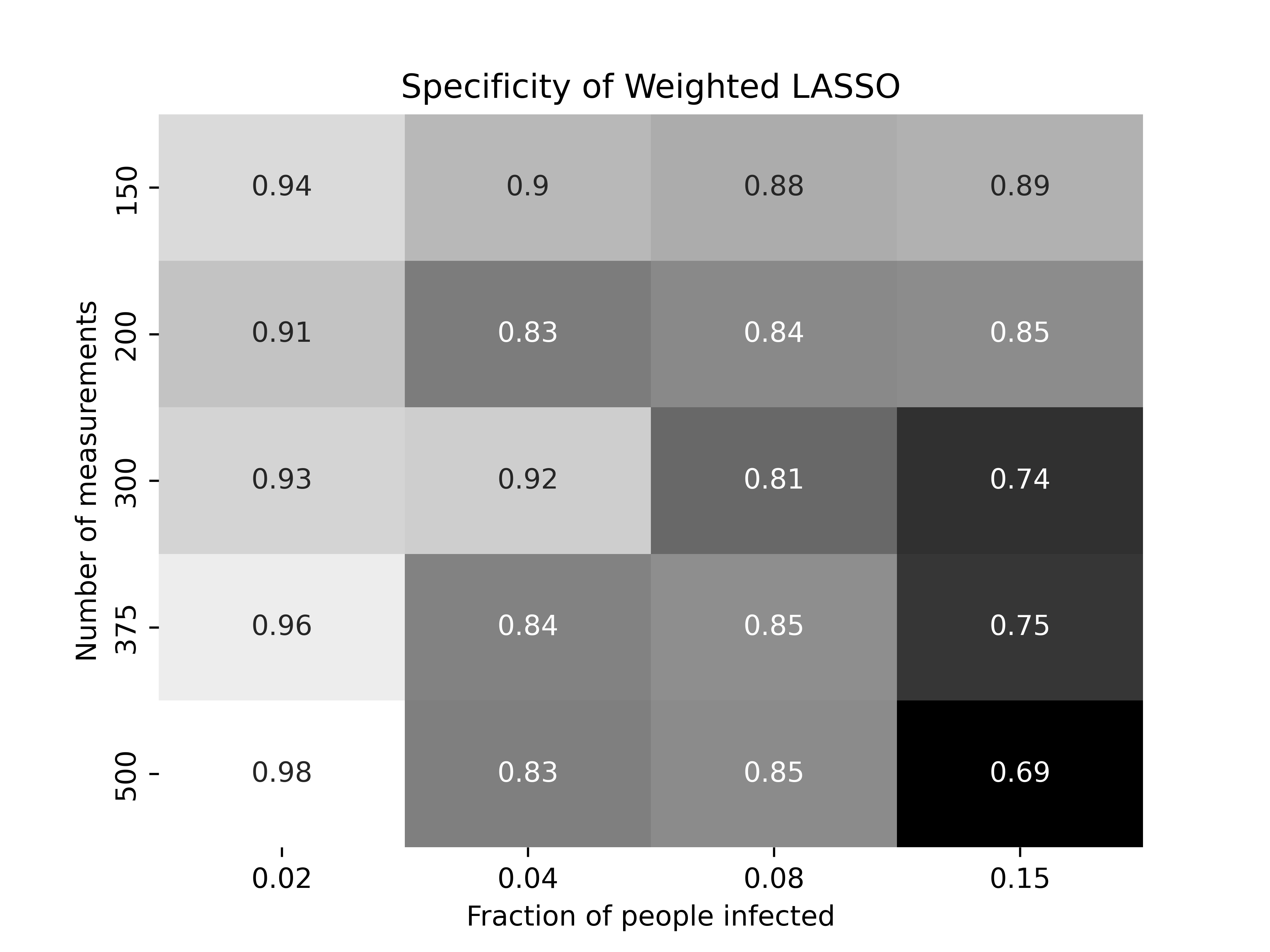}
  \label{fig:wlasso_spec_0.5}
\end{subfigure}
\vspace{-1\baselineskip}
\caption{RRMSE, Sensitivity and Specificity values for $q=0.5$, versus the number of measurements $n$ and sparsity (fraction $f_s$ of people infected). The top and bottom rows represent results for \textsc{Lasso} and \textsc{WLasso}, respectively. Higher is better for Sensitivity and Specificity, and lower is better for RRMSE. Additional results with $\sigma \in \{0.1,0.15,0.2\}$ can be found in Sec.4 of the Supplemental Material.}
\label{fig:mainres}
\end{figure}

\begin{figure}[t!]
    \centering 
\begin{subfigure}{0.25\textwidth}
  \includegraphics[width=1.1\linewidth]{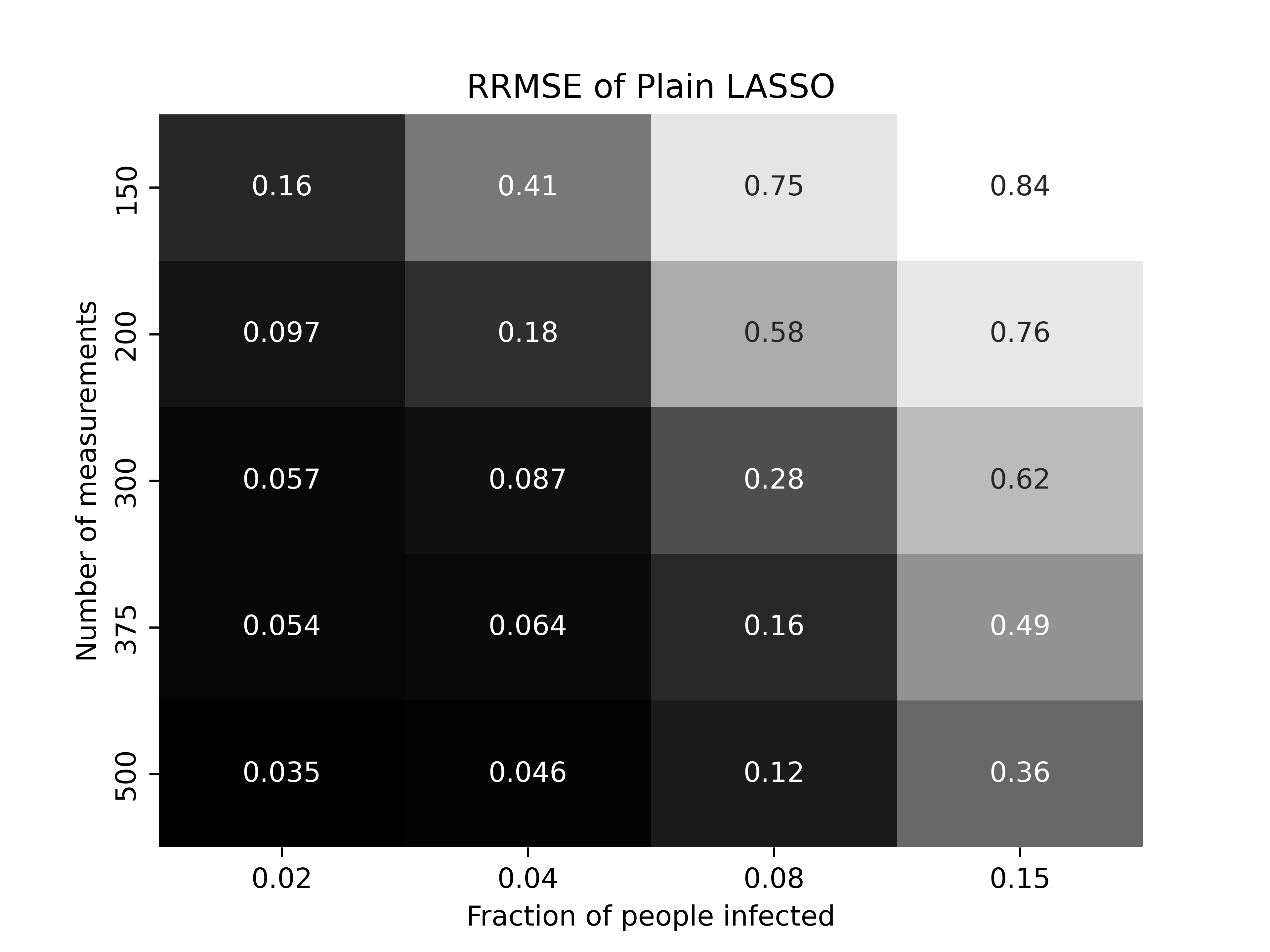}
  \caption{$q=0.05$}
  \label{fig:1}
\end{subfigure}\hfil 
\begin{subfigure}{0.25\textwidth}
  \includegraphics[width=1.1\linewidth]{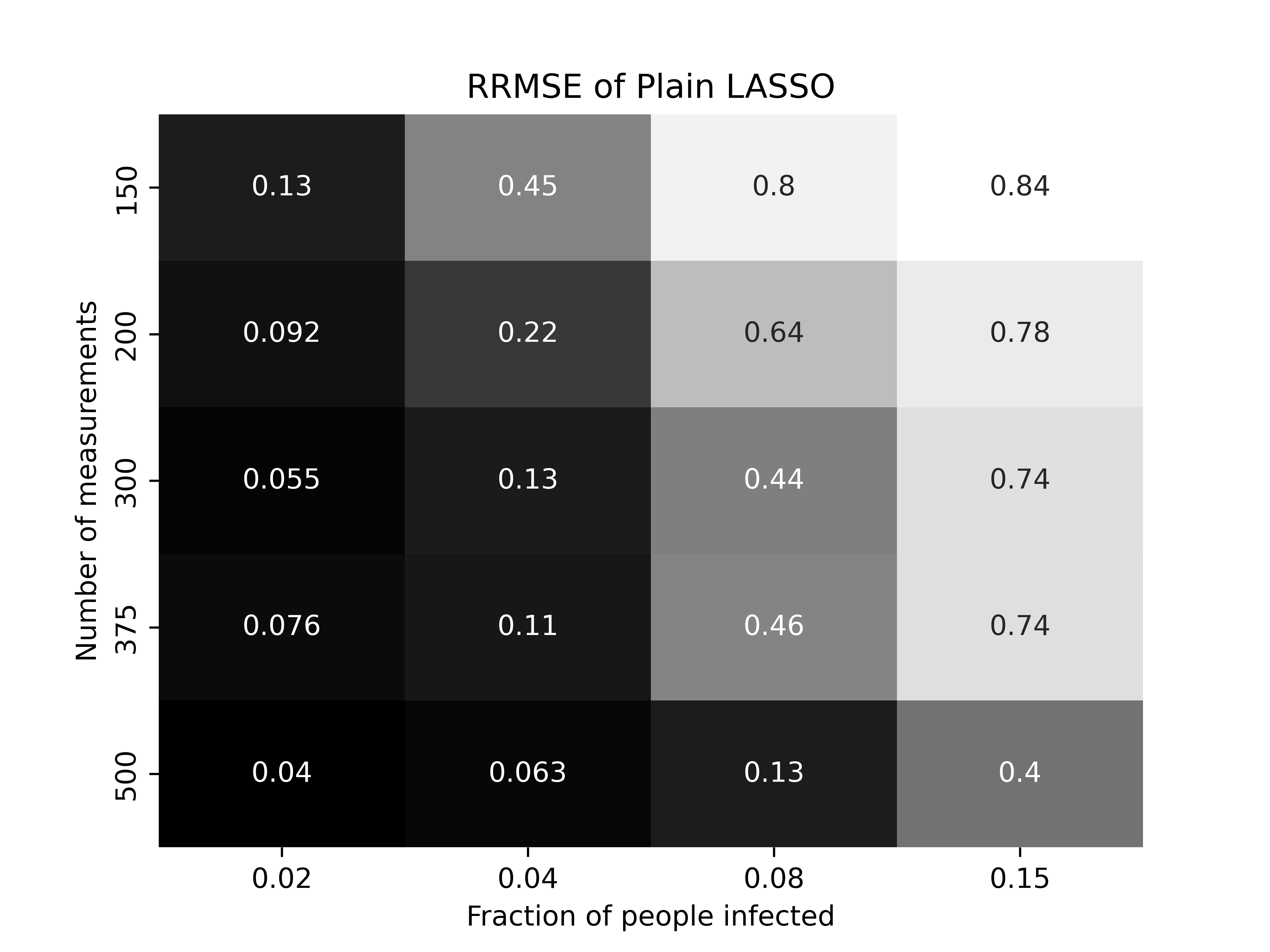}
  \caption{$q=0.1$}
  \label{fig:2}
\end{subfigure}\hfil 
\begin{subfigure}{0.25\textwidth}
  \includegraphics[width=1.1\linewidth]{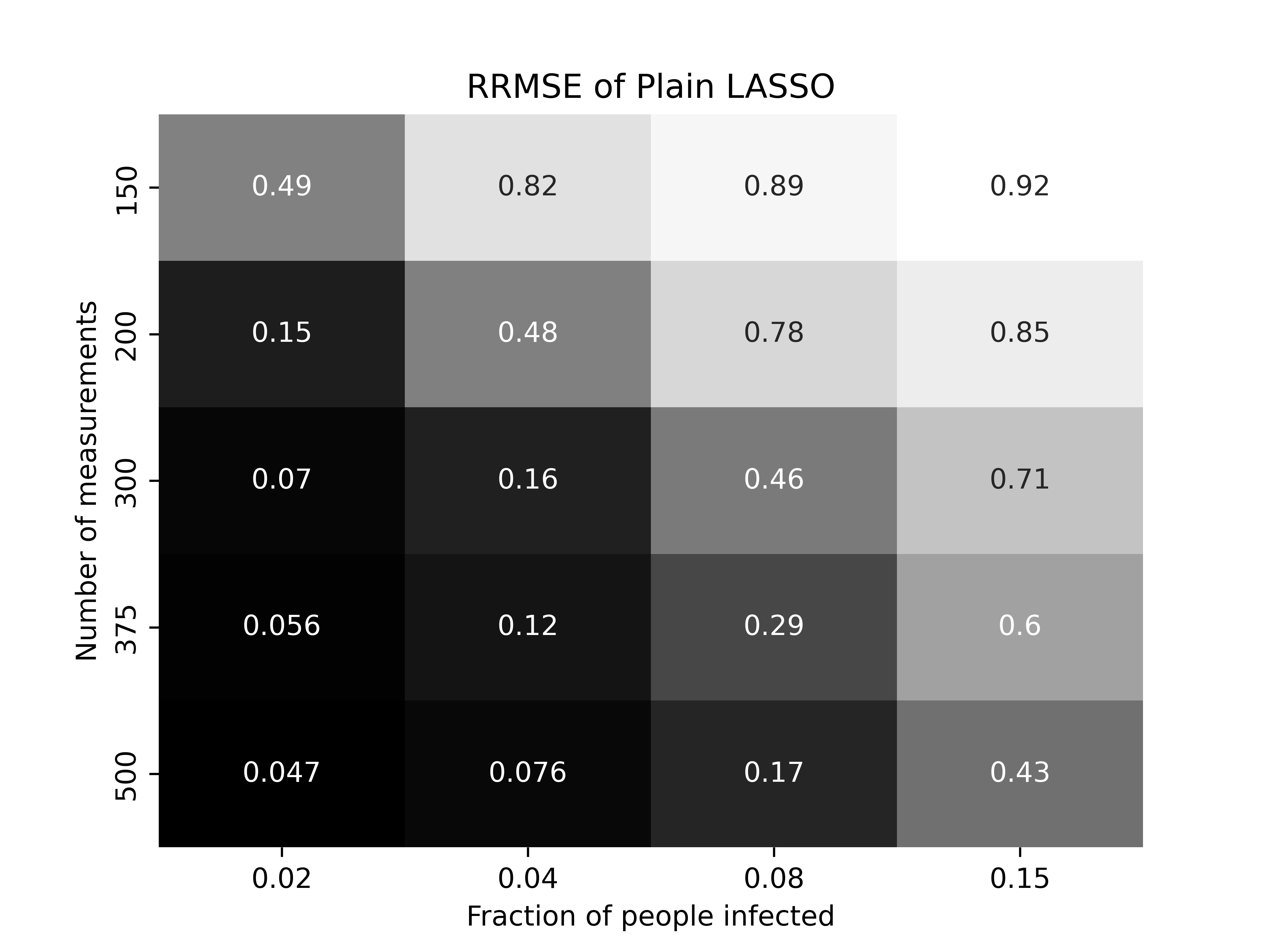}
  \caption{$q=0.25$}
  \label{fig:3}
\end{subfigure}\hfil 
\begin{subfigure}{0.25\textwidth}
  \includegraphics[width=1.1\linewidth]{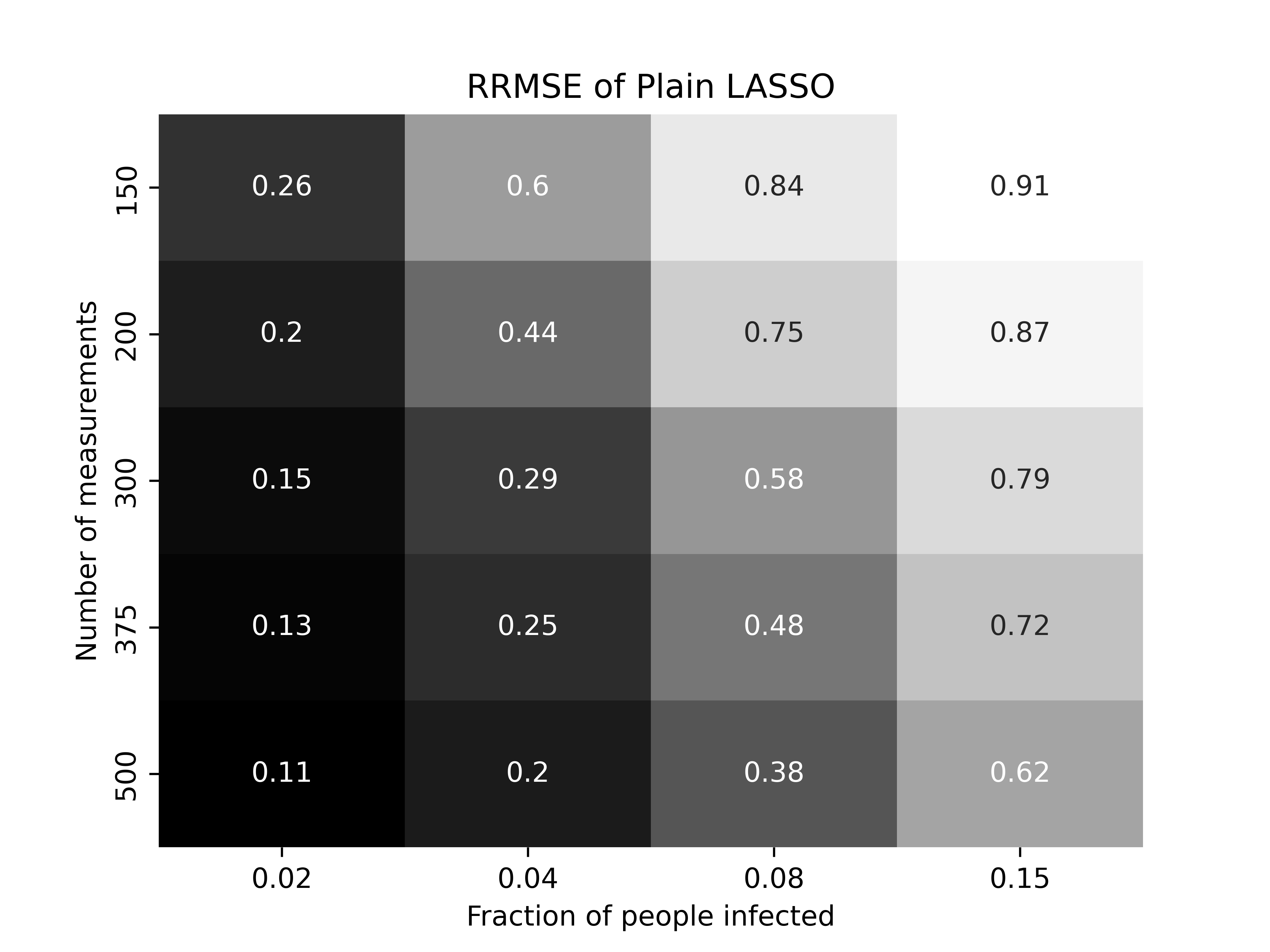}
  \caption{$q=0.75$}
  \label{fig:4}
\end{subfigure}

\begin{subfigure}{0.25\textwidth}
  \includegraphics[width=1.1\linewidth]{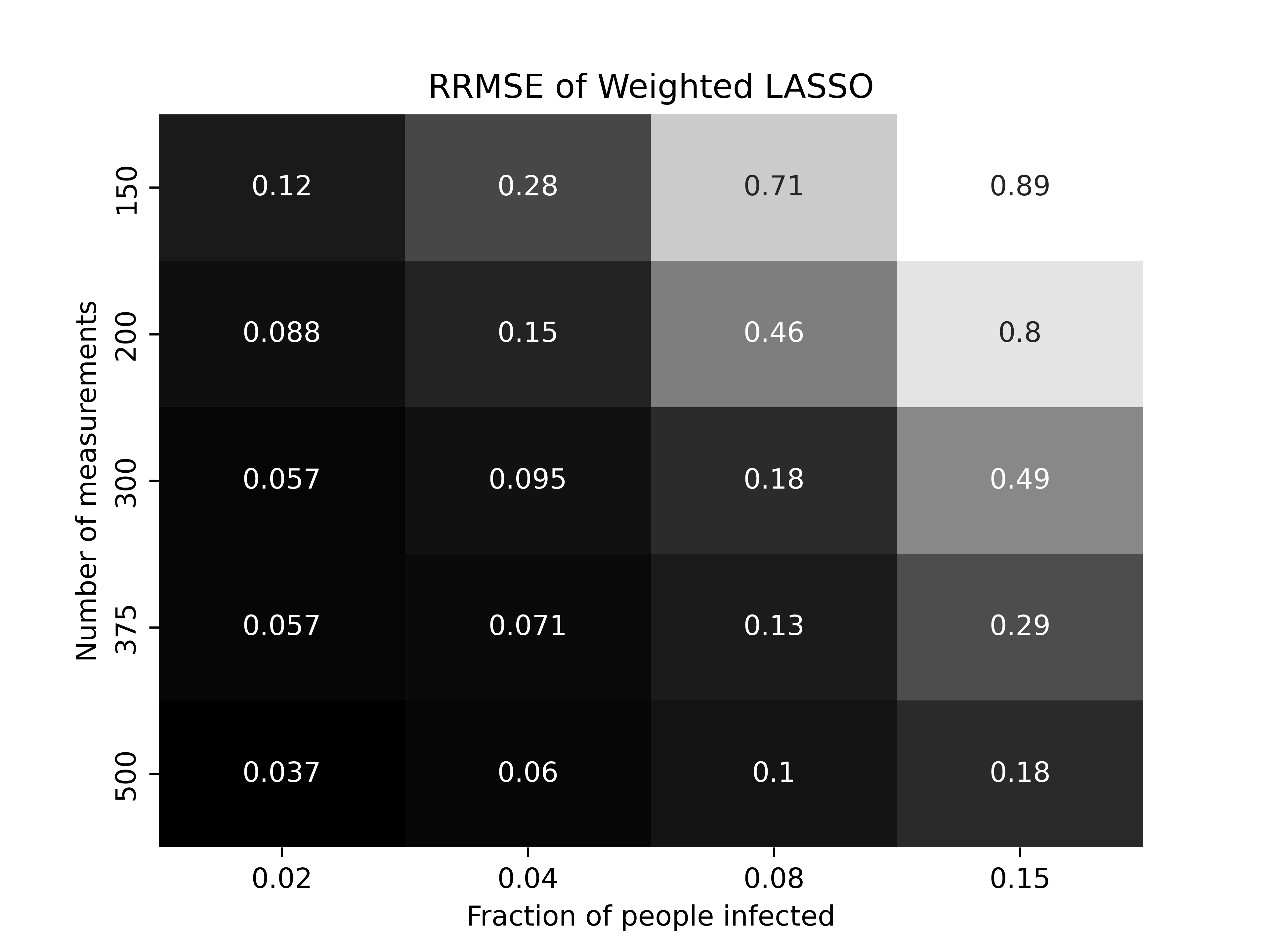}
  \caption{$q=0.05$}
  \label{fig:5}
\end{subfigure}\hfil
\begin{subfigure}{0.25\textwidth}
  \includegraphics[width=1.1\linewidth]{newplots/perf-100-0.95-0.1-rrmse.png}
  \caption{$q=0.1$}
  \label{fig:6}
\end{subfigure}\hfil
\begin{subfigure}{0.25\textwidth}
  \includegraphics[width=1.1\linewidth]{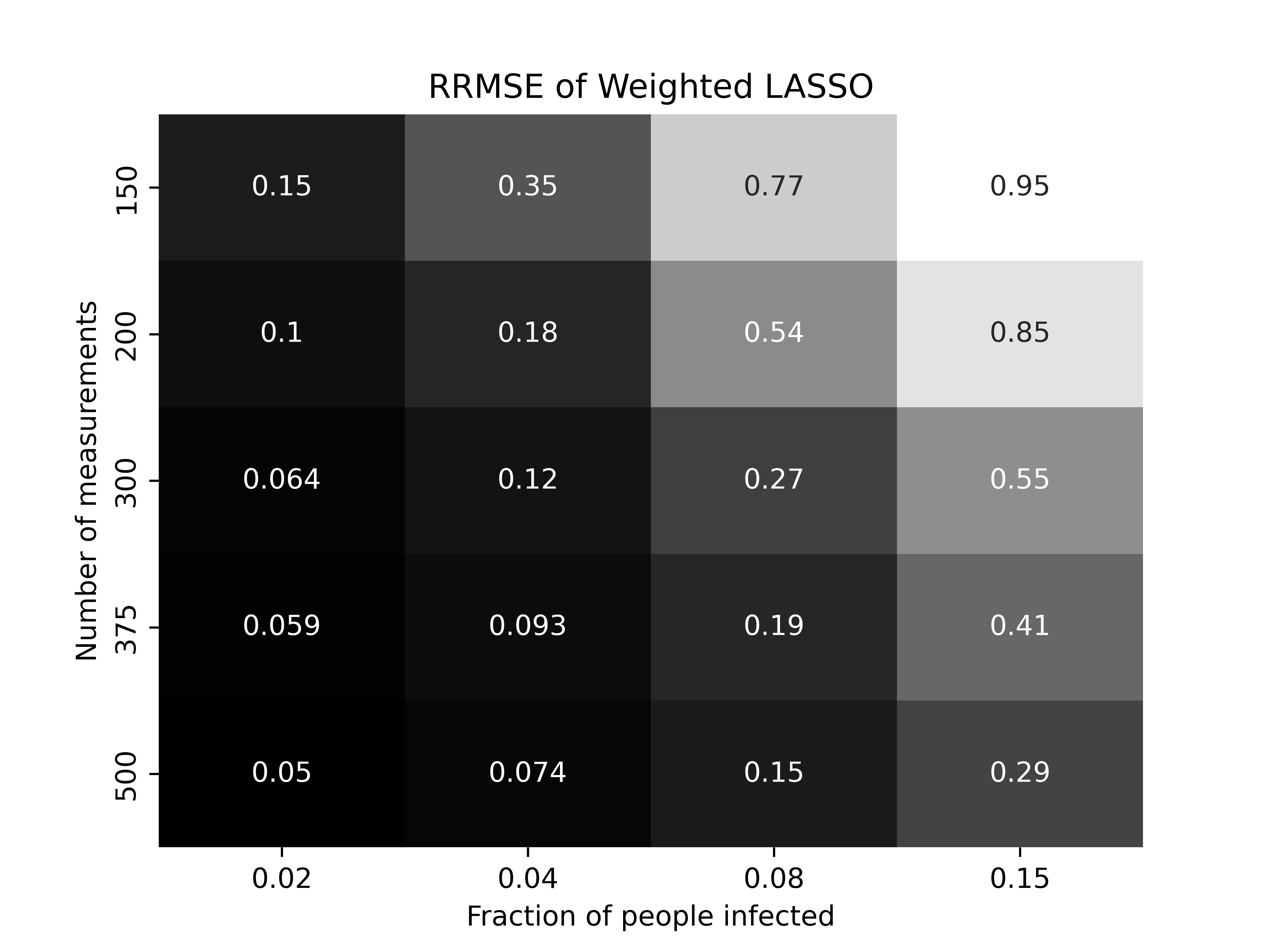}
  \caption{$q=0.25$}
  \label{fig:7}
\end{subfigure}\hfil
\begin{subfigure}{0.25\textwidth}
  \includegraphics[width=1.1\linewidth]{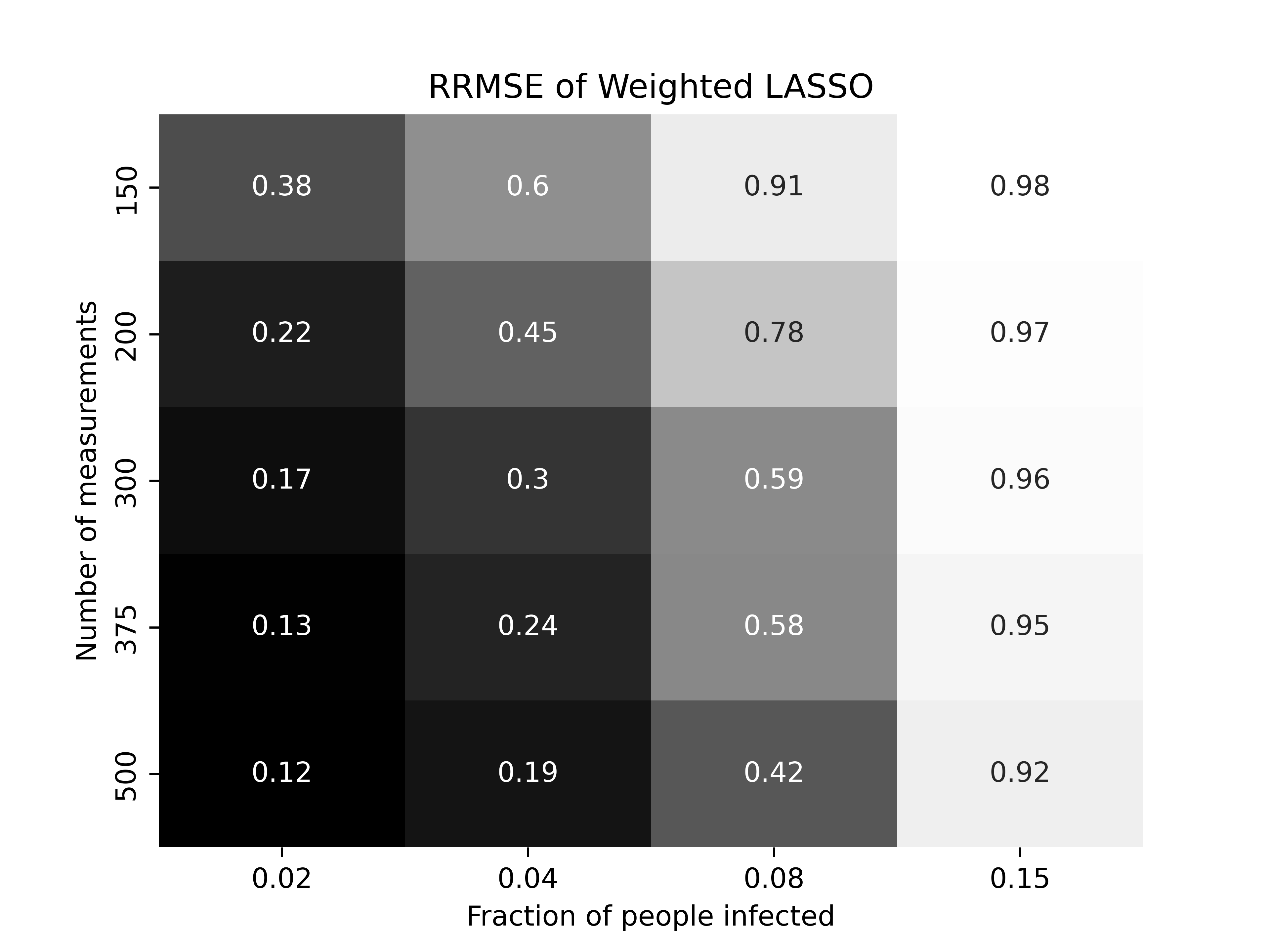}
  \caption{$q=0.75$}
  \label{fig:8}
\end{subfigure}
\caption{RRMSE values (lower is better) for $q \in \{0.05,0.1,0.25,0.75\}$ (left to right), versus the number of measurements $n$ and $f_s$ (fraction of people infected out of a total of $p$). The top and bottom rows represent the results for \textsc{Lasso} and \textsc{WLasso}, respectively.}
\label{fig:rrmse_diff_q}
\end{figure}

We evaluate the derived expressions for the weights used in the \textsc{Lasso} and \textsc{WLasso} estimators, using simulated signal reconstruction in pooled RT-PCR. In particular, we use the \texttt{cvxpy} \cite{cvxpy1,cvxpy2} Python package to solve the \textsc{WLasso} problem using inputs $\bm{\Tilde{y}}, \bm{\Tilde{A}}$, which are derived from $\bm{y}, \bm{A}$ as explained in Sec.~\ref{ssec:rescalerecenter}. In our simulations, we take the standard deviation of the (multiplicative) Gaussian noise $\sigma$ to be 0.05, which is in tune with assumption A3 of Proposition 1 (see Comment 1 after Propositions 1 and 2). Additional results with $\sigma \in \{0.1,0.15,0.2\}$ can be found in Sec.4 of the Supplemental Material. The amplification factor $q_a$ is taken to be $0.95$ (a reasonable choice since this factor is known to be close to 1 as the molecules roughly double in each RT-PCR cycle). Note that the values $\sigma, q, q_a$ need not be specified as input to \textsc{Lasso} or \textsc{WLasso} and are needed only for problem specification. The value of $p$ (the number of samples being tested) is fixed at 1000 and the number of pools $n$ is varied among $\{150,200,300,375,500\}$. The pooling matrix $\bm{A}$ is a random Bernoulli matrix with each entry equal to $1$ with probability $q$, and $0$ otherwise, freshly sampled for each variation. We experiment with values of $q$ in $\{0.05, 0.1, 0.25, 0.5, 0.75\}$. We mostly focus on the case of $q=0.5$ here - the complete set of plots will be made publicly available on GitHub along with the code used for simulation\footnote{\url{https://github.com/sudoRicheek/Weighted-LASSO-MGN-Simulations}}. The underlying ground truth signal $\bm{x^*}$ is chosen to have an $\ell_0$ norm given by $f_s \times p$ where the fraction $f_s$ is sampled from $\{0.02, 0.04, 0.08, 0.15\}$. The significant non-zero entries of $\bm{x^*}$ are uniformly randomly sampled from the interval $[1,1000]$, which reasonably represents the viral loads of infected individuals. 
The indices of the significant non-zero entries are drawn independently from a uniform distribution.
The other entries are uniformly randomly sampled from $[0,0.2]$. A viral load of 0.2 is used as a threshold to determine infected individuals from the reconstructed viral loads. However, slight variations in this threshold do not alter the Sensitivity and Specificity values stated in \cref{fig:mainres}.

\begin{figure}[t!]
    \centering 
\begin{subfigure}{0.33\textwidth}
  \includegraphics[width=1.1\linewidth]{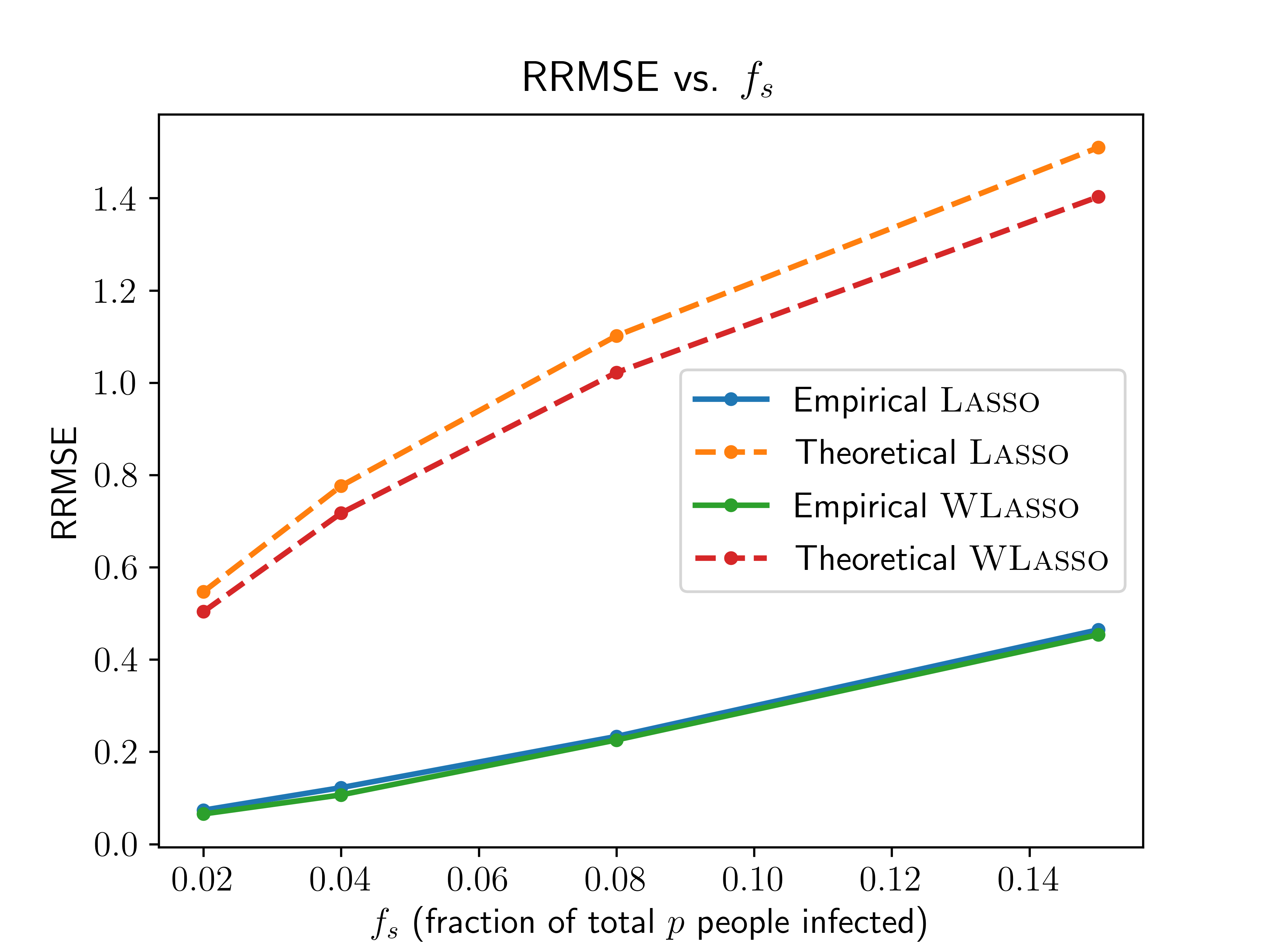}
  \label{fig:trends}
\end{subfigure}\hfil 
\begin{subfigure}{0.33\textwidth}
  \includegraphics[width=1.1\linewidth]{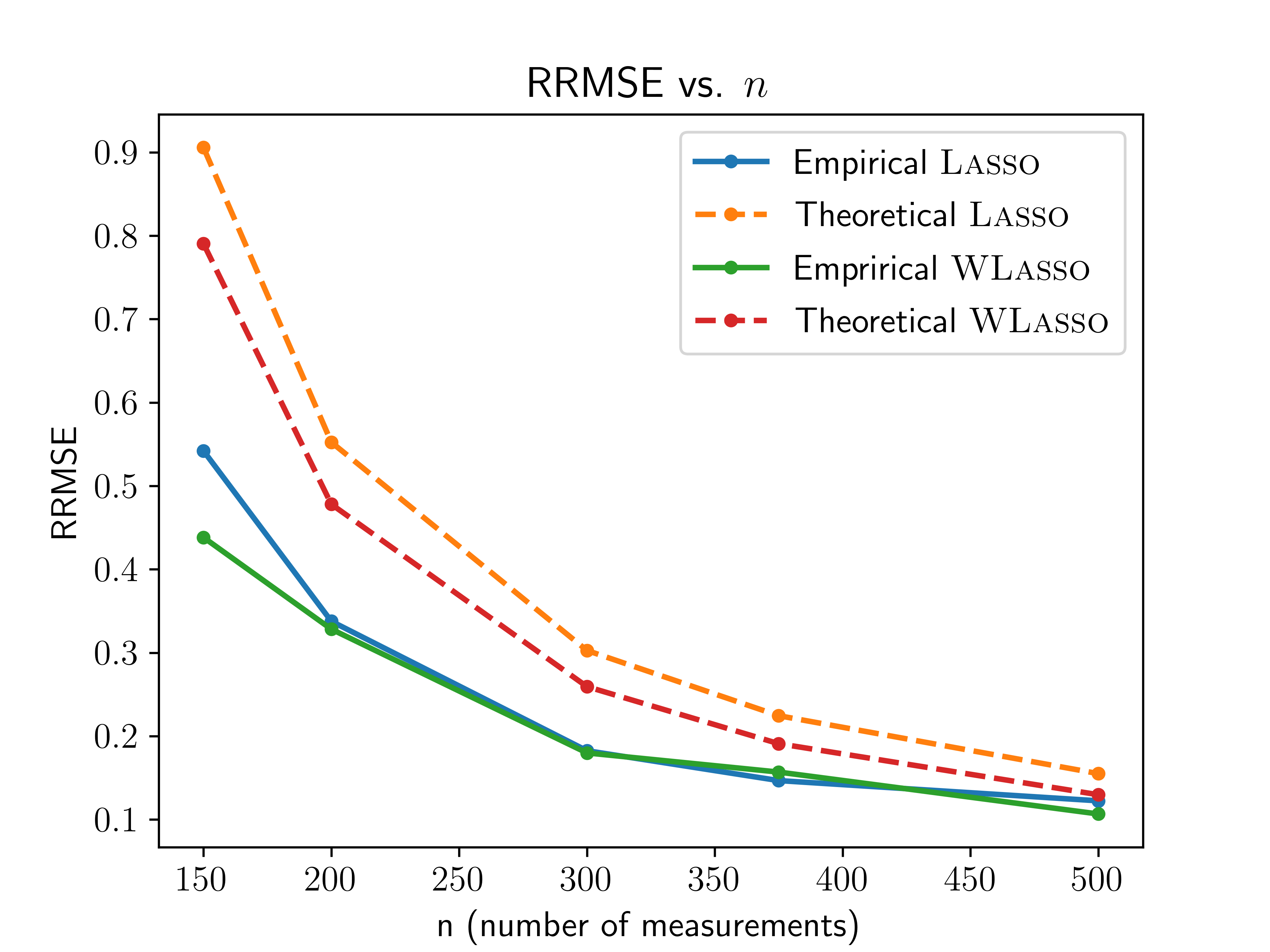}
  \label{fig:trendn}
\end{subfigure}\hfil 
\begin{subfigure}{0.33\textwidth}
  \includegraphics[width=1.1\linewidth]{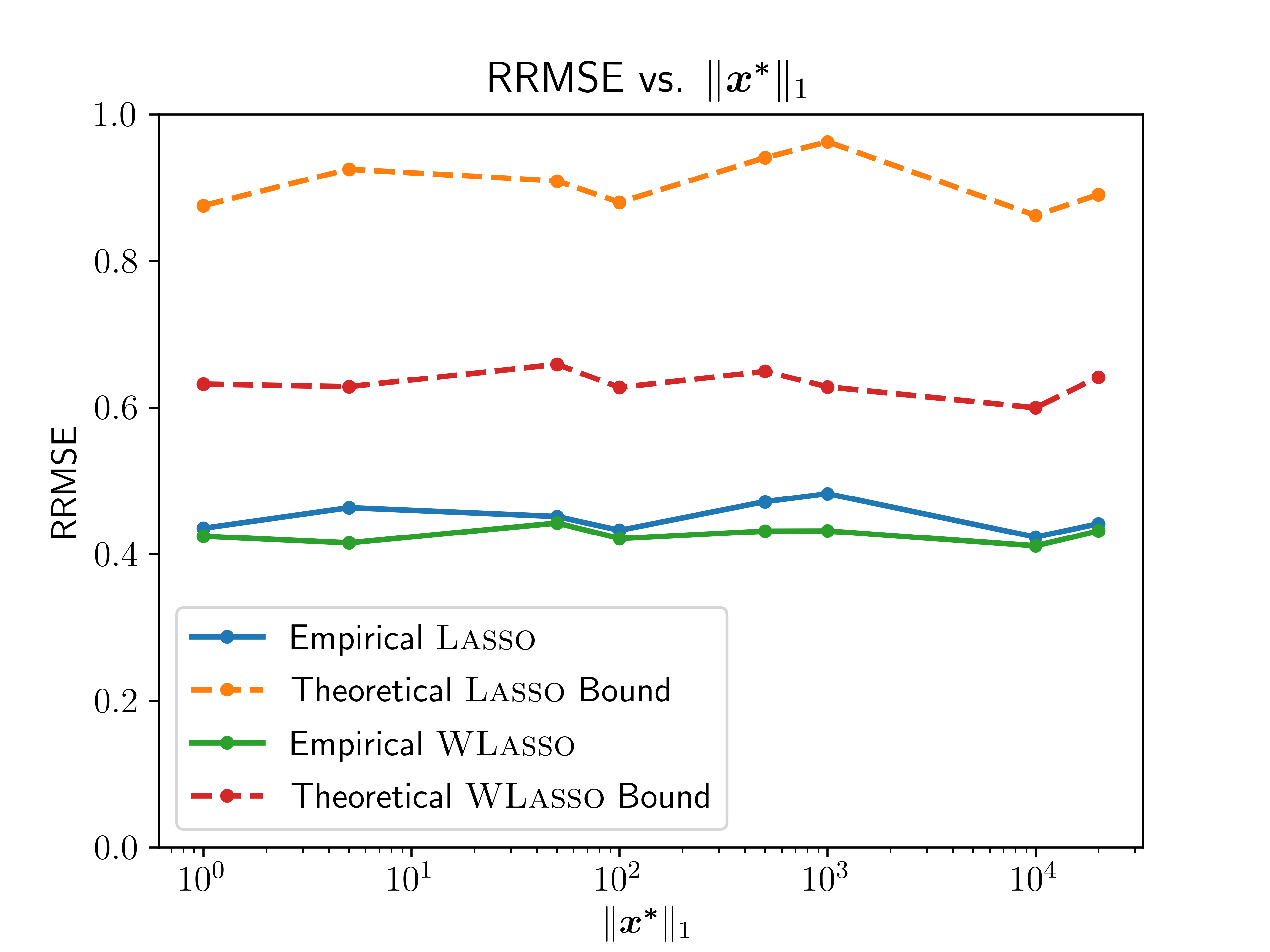}
  \label{fig:trendx1}
\end{subfigure}
\vspace{-1\baselineskip}
\caption{General trends of RRMSE vs $\{f_s, n, \|\bm{x^*}\|_1\}$ (left to right). In each case, one of the parameters is varied, keeping the rest fixed to default values.}
\label{fig:trends}
\end{figure}

For each combination of $n, f_s, q$, we execute 200 independent simulations and calculate averages of three intuitive and widely used measures: the RRMSE (Relative Root Mean-Squared-Error), the Sensitivity, and the Specificity, across the runs. The RRMSE is computed as $\dfrac{\|\bm{x^*}-\bm{\hat{x}}\|_2}{\|\bm{x^*}\|_2}$, where $\bm{\hat{x}}$ represents the estimate of the ground truth vector $\bm{x}$. Let TP, TN, FP, FN denote the number of true positives, true negatives, false positives and false negatives respectively. The sensitivity is defined as $\dfrac{TP}{TP+FN}$ (the probability of a positive test result for an unhealthy sample), whereas the specificity is defined as $\dfrac{TN}{TN+FP}$ (the probability of a negative test result for a healthy sample). Higher values are better for sensitivity and specificity, whereas the opposite holds true for RRMSE. We report the measures as grid maps for both  \textsc{Lasso} and \textsc{WLasso}, with $n$ and $f_s$ being the independent variables selected from the ranges mentioned earlier. For each combination, we calculate a regularization factor $\gamma$ to be multiplied to the weight $\beta$ for \textsc{Lasso} or the weights $\{\beta_k\}_{k=1}^p$ for \textsc{WLasso}, using cross-validation on 5 preliminary runs.

We show the plots for the RRMSE values $(q=0.5)$ in the first column of \cref{fig:mainres}. As expected, the RRMSE increases with the $f_s$ (as a greater amount of information is now compressed down to the same dimensions) and decreases with the number of measurements (as more measurements capture more information). Out of \textsc{Lasso} and \textsc{Wlasso}, neither estimator is consistently superior to the other. This is, however unsurprising as the bounds we derive on the weights are worst-case bounds and may not apply to the average case. Nevertheless, the values obtained clearly demonstrate the suitability of the derived weights for reconstruction. We show the corresponding values for the Sensitivity and Specificity (for $q=0.5$) in the middle and rightmost columns of \cref{fig:mainres} respectively. Interestingly, for higher fraction of people infected $(f_s)$, \textsc{WLasso} achieves a better specificity while being competitive with \textsc{Lasso} in the sensitivity values. Both methods achieve high values on both measures at lower sparsity levels.

\begin{figure}[t!]
    \centering
    \includegraphics[width=0.7\linewidth]{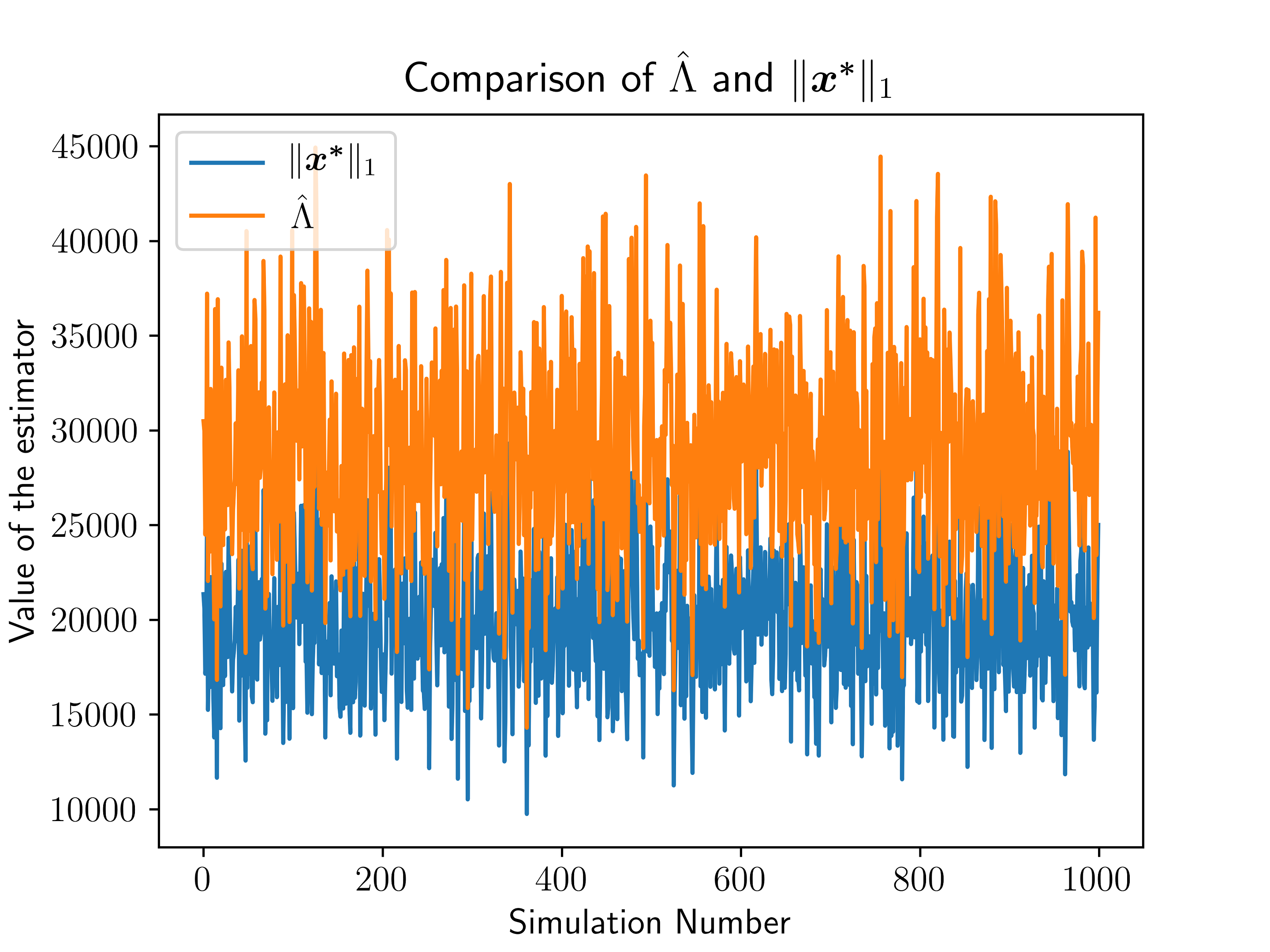}
    \caption{Comparison of the estimator $\hat{\Lambda}$ with $\|\bm{x^*}\|_1$ over a set of 1000 simulation runs (corresponding values of the two quantities are plotted together).}
    \label{fig:lambdahatcomp}
\end{figure}

We also compare the obtained results with the cases when the value of $q$ for the sensing matrix is chosen to be from $\{0.05, 0.1, 0.25, 0.75\}$. These results are shown in \cref{fig:rrmse_diff_q}. As can be seen, smaller values of $q$ yield lower RRMSE values as the average number of participants per pool decreases with $q$. These trends will hold true as long as every sample participates in two or more pools and no pool is left empty. 

As can be seen from Propositions \ref{prop1}, \ref{prop2} and \ref{prop3}, the RRMSE value scale as $O(\sqrt{s})$ w.r.t. $s$, as $O(1/n)$ w.r.t. $n$ and as $O(1)$ w.r.t. $\|\bm{x^*}\|_1$, keeping all other factors constant and set to their default values. We verified experimentally that these trends hold true, ignoring constant factors. The verification was performed over 500 independently generated signals. These trends are demonstrated in \cref{fig:trends}. 

Finally, the results in Propositions \ref{prop1}, \ref{prop2} and \ref{prop3} use an estimate of $\|\bm{x^*}\|_1$, denoted as $\hat{\Lambda}$. In \cref{fig:lambdahatcomp}, we demonstrate that $\hat{\Lambda}$ is a good estimate for $\|\bm{x^*}\|_1$ and its value empirically never crosses $2\|\bm{x^*}\|_1$. This is shown for 1000 independent simulations. As an example, the average value of $\|\bm{x^*}\|_1$ across 1000 simulations was $\sim 19,900$ with a standard deviation of $\sim 3500$, whereas the 
average value of $\hat{\Lambda}$ was $\sim 28,900$ with a standard deviation of $\sim 5200$. 

\noindent\textbf{Discussion:} In our group's previous work \cite{Ghosh2021}, we have already shown comparisons between many different estimators for compressive RT-PCR. These estimators include \textsc{Lasso}, sparse Bayesian learning (\textsc{Sbl}), orthogonal matching pursuit (\textsc{Omp}) among others. Out of these, \textsc{Lasso} is very simple to implement and performs empirically very well. However, there are no performance bounds for \textsc{Lasso} given this noise model (multiplicative log-normal) in RT-PCR in the literature. The aim of our paper is to fill that gap. Apart from this, there exists very little literature on compressive inversion given multiplicative noise in the measurements, \cite{Zhou2022} being a recent reference. However, the noise model in \cite{Zhou2022} is multiplicative \emph{Gaussian} (which may be positive or negative), whereas the noise model in our paper is multiplicative \emph{lognormal} (which is always non-negative). Hence, no comparison is possible in such a case. There certainly exist many algorithms for image filtering in multiplicative noise, however our work is on compressive inversion, and hence direct comparisons are not possible. Further note that first filtering $\boldsymbol{y}$ and then reconstructing $\boldsymbol{x^*}$ will not work well, as $\boldsymbol{y}$ is not a sparse vector unlike $\boldsymbol{x^*}$ which would lead to poor filtering performance.

\section{Conclusion}
Despite extensive work on the theoretical guarantees for compressed sensing in the additive Gaussian Noise setting, current literature lacks performance bounds for the case of multiplicative lognormal Noise which describes the noise in the scenario of RT-PCR Group Testing. In this work we develop performance bounds for the weighted and unweighted \textsc{Lasso} problems in the case of multiplicative lognormal noise, and derive expressions for weights corresponding to these bounds. We then show via simulation that our weights yield good performance in practice. We hope that our efforts inspire future work on more varied as well as general noise models in compressed sensing. A full analysis of performance bounds for the noise model without use of Taylor expansion, along the lines of \cite{Rishav2021}, is an important avenue for future work. Our current work considers randomly generated Bernoulli matrices as pooling matrices. Extending our analysis to the important case of deterministic matrices \cite{Ghosh2021} is another useful direction for future work. Multiplicative noise, albeit of a different form, also shows up in some other modalities such as synthetic aperture radar \cite{Pournaghshband2023,Aubert2008,Zhou2022} and photo-acoustic tomography \cite[Equation 17]{Arridge2016}. Developing compressed sensing bounds in such modalities is also another future research direction. 








\bibliographystyle{elsarticle-num}
\bibliography{bibtex}

\begin{thebibliography}{10}
\expandafter\ifx\csname url\endcsname\relax
  \def\url#1{\texttt{#1}}\fi
\expandafter\ifx\csname urlprefix\endcsname\relax\def\urlprefix{URL }\fi
\expandafter\ifx\csname href\endcsname\relax
  \def\href#1#2{#2} \def\path#1{#1}\fi

\bibitem{Ghosh2021}
S.~Ghosh, R.~Agarwal, M.~A. Rehan, S.~Pathak, P.~Agarwal, Y.~Gupta, S.~Consul,
  N.~Gupta, Ritika, R.~Goenka, A.~Rajwade, M.~Gopalkrishnan, A compressed
  sensing approach to pooled {RT}-{PCR} testing for {COVID}-19 detection,
  {IEEE} Open Journal of Signal Processing 2 (2021) 248--264.
\newblock \href {https://doi.org/10.1109/ojsp.2021.3075913}
  {\path{doi:10.1109/ojsp.2021.3075913}}.

\bibitem{Hunt2018}
X.~J. Hunt, P.~Reynaud-Bouret, V.~Rivoirard, L.~Sansonnet, R.~Willett, A
  data-dependent weighted {LASSO} under {P}oisson noise, IEEE Transactions on
  Information Theory 65~(3) (2018) 1589--1613.
\newblock \href {https://doi.org/10.1109/TIT.2018.2869578}
  {\path{doi:10.1109/TIT.2018.2869578}}.

\bibitem{Jawerth2020}
N.~Jawerth, How is the {COVID-19} virus detected using real time {RT-PCR}?,
  \url{https://www.iaea.org/newscenter/news/how-is-the-covid-19-virus-detected-using-real-time-rt-pcr}.

\bibitem{PooledWiki2021}
List of countries implementing pool testing strategy against {COVID}-19,
  \url{https://en.wikipedia.org/wiki/List_of_countries_implementing_pool_testing_strategy_against_COVID-19},
  last retrieved, Oct 2021.

\bibitem{Shental2020}
N.~Shental, et~al., Efficient high throughput {SARS-CoV-2} testing to detect
  asymptomatic carriers, Sci. Adv. 6~(37) (Sep. 2020).
\newblock \href {https://doi.org/10.1126/sciadv.abc5961}
  {\path{doi:10.1126/sciadv.abc5961}}.

\bibitem{Heiderzadeh2021}
A.~Heidarzadeh, K.~Narayanan, Two-stage adaptive pooling with {RT-qPCR} for
  {COVID-19} screening, in: ICASSP, 2021.

\bibitem{Goenka2021}
R.~Goenka, S.~Cao, C.~Wong, A.~Rajwade, D.~Baron,
  \href{https://arxiv.org/abs/2106.02699}{Contact tracing information improves
  the performance of group testing algorithms} (2021).
\newline\urlprefix\url{https://arxiv.org/abs/2106.02699}

\bibitem{Goenka2021contact}
R.~Goenka, S.-J. Cao, C.-W. Wong, A.~Rajwade, D.~Baron, Contact tracing
  enhances the efficiency of covid-19 group testing, in: ICASSP 2021-2021 IEEE
  International Conference on Acoustics, Speech and Signal Processing (ICASSP),
  IEEE, 2021, pp. 8168--8172.

\bibitem{Dorfman1943}
R.~Dorfman, {The Detection of Defective Members of Large Populations}, The
  Annals of Mathematical Statistics 14~(4) (1943) 436 -- 440.
\newblock \href {https://doi.org/10.1214/aoms/1177731363}
  {\path{doi:10.1214/aoms/1177731363}}.

\bibitem{israeli_news}
\href{https://www.israel21c.org/israelis-introduce-method-for-accelerated-covid-19-testing/}{Israelis
  introduce a method for accelerated covid-19 testing}.
\newline\urlprefix\url{https://www.israel21c.org/israelis-introduce-method-for-accelerated-covid-19-testing/}

\bibitem{healthcare_europe}
\href{https://healthcare-in-europe.com/en/news/corona-pool-testing-increases-worldwide-capacities-many-times-over.html}{Corona
  pool testing increases worldwide capacities many times over}.
\newline\urlprefix\url{https://healthcare-in-europe.com/en/news/corona-pool-testing-increases-worldwide-capacities-many-times-over.html}

\bibitem{Larremore2021}
D.~B. Larremore, B.~Wilder, E.~Lester, S.~Shehata, J.~M. Burke, J.~A. Hay,
  M.~Tambe, M.~Mina, R.~Parker, Test sensitivity is secondary to frequency and
  turnaround time for {COVID-19} screening, Science Advances 7~(1) (2021).

\bibitem{Donoho2006}
D.~Donoho, Compressed sensing, IEEE Transactions on Information Theory 52~(4)
  (2006) 1289--1306.
\newblock \href {https://doi.org/10.1109/TIT.2006.871582}
  {\path{doi:10.1109/TIT.2006.871582}}.

\bibitem{Aldridge2019}
M.~Aldridge, O.~Johnson, J.~Scarlett, et~al., Group testing: an information
  theory perspective, Foundations and Trends{\textregistered} in Communications
  and Information Theory 15~(3-4) (2019) 196--392.

\bibitem{THW2015}
T.~Hastie, R.~Tibshirani, M.~Wainwright, Statistical Learning with Sparsity:
  The {LASSO} and Generalizations, {CRC} Press, 2015.

\bibitem{Huang2010}
J.~Huang, T.~Zhang, {The benefit of group sparsity}, The Annals of Statistics
  38~(4) (2010) 1978 -- 2004.

\bibitem{Li2018}
Y.~Li, G.~Raskutti, Minimax optimal convex methods for {P}oisson inverse
  problems under $\ell_{q}$ -ball sparsity, IEEE Transactions on Information
  Theory 64~(8) (2018) 5498--5512.
\newblock \href {https://doi.org/10.1109/TIT.2018.2850365}
  {\path{doi:10.1109/TIT.2018.2850365}}.

\bibitem{Negahban2012}
S.~Negahban, P.~Ravikumar, M.~J. Wainwright, B.~Yu, A unified framework for
  high-dimensional analysis of {M}-estimators with decomposable regularizers,
  Statistical Science 27~(4) (2012) 538--557.

\bibitem{Rishav2021}
R.~Chatterjee, Reconstruction of sparse signals under the influence of
  multiplicative log-normal noise using l1-regularization,
  \url{https://tinyurl.com/4w8rr5mz}, {MSc} Thesis, IIT Bombay (2021).

\bibitem{Rajwade2021}
A.~Rajwade, K.~S. Gurumoorthy, Two penalized estimators based on variance
  stabilization transforms for sparse compressive recovery with {P}oisson
  measurement noise, Signal Processing 188 (2021) 108186.

\bibitem{Bohra2019}
P.~Bohra, D.~Garg, K.~S. Gurumoorthy, A.~Rajwade, Variance-stabilization-based
  compressive inversion under {P}oisson or {P}oisson--{G}aussian noise with
  analytical bounds, Inverse Problems 35~(10) (2019) 105006.

\bibitem{Stojanovic2016}
V.~Stojanovic, N.~Nedic, D.~Prsic, L.~Dubonjic, Optimal experiment design for
  identification of arx models with constrained output in non-gaussian noise,
  Applied Mathematical Modelling 40~(13-14) (2016) 6676--6689.

\bibitem{Stojanovic2016_2}
V.~Stojanovic, N.~Nedic, Robust identification of oe model with constrained
  output using optimal input design, Journal of the Franklin Institute 353~(2)
  (2016) 576--593.

\bibitem{surveygt}
E.~A. Daniel, et~al., {{P}ooled {T}esting {S}trategies for
  {S}{A}{R}{S}-{C}o{V}-2 diagnosis: {A} comprehensive review}, Diagn Microbiol
  Infect Dis 101~(2) (2021) 115432.

\bibitem{Karlen2007}
Y.~Karlen, A.~McNair, S.~Perseguers, C.~Mazza, N.~Mermod, Statistical
  significance of quantitative {PCR}, BMC bioinformatics 8 (2007) 1--16.

\bibitem{Arildsen2014}
T.~Arildsen, T.~Larsen, Compressed sensing with linear correlation between
  signal and measurement noise, Signal Processing 98 (2014) 275--283.

\bibitem{Pollard1979}
J.~H. Pollard, A handbook of numerical and statistical techniques: with
  examples mainly from the life sciences, CUP Archive, 1979.

\bibitem{Gopal2020}
P.~Gopal, S.~Chandran, I.~Svalbe, A.~Rajwade, Low radiation tomographic
  reconstruction with and without template information, Signal Processing 175
  (2020) 107582.

\bibitem{Dupe2009}
F.-X. Dup{\'e}, J.~M. Fadili, J.-L. Starck, A proximal iteration for
  deconvolving {P}oisson noisy images using sparse representations, IEEE
  Transactions on Image Processing 18~(2) (2009) 310--321.

\bibitem{BernoulliInequality_wiki}
Bernoulli's inequality,
  \url{https://en.wikipedia.org/wiki/Bernoulli's_inequality}.

\bibitem{Brannan2006}
D.~A. Brannan, A first course in mathematical analysis, Cambridge University
  Press, 2006.

\bibitem{cvxpy1}
S.~Diamond, S.~Boyd, {CVXPY}: {A} {P}ython-embedded modeling language for
  convex optimization, Journal of Machine Learning Research 17~(83) (2016)
  1--5.

\bibitem{cvxpy2}
A.~Agrawal, R.~Verschueren, S.~Diamond, S.~Boyd, A rewriting system for convex
  optimization problems, Journal of Control and Decision 5~(1) (2018) 42--60.

\bibitem{Zhou2022}
W.~Zhou, S.~Jalali, A.~Maleki, Compressed sensing in the presence of speckle
  noise, IEEE Transactions on Information Theory 68~(10) (2022) 6964--6980.

\bibitem{Pournaghshband2023}
R.~Pournaghshband, M.~Modarres-Hashemi, A novel block compressive sensing
  algorithm for {SAR} image formation, Signal Processing 210 (2023) 109053.

\bibitem{Aubert2008}
G.~Aubert, J.-F. Aujol, A variational approach to removing multiplicative
  noise, SIAM journal on applied mathematics 68~(4) (2008) 925--946.

\bibitem{Arridge2016}
S.~Arridge, P.~Beard, M.~Betcke, B.~Cox, N.~Huynh, F.~Lucka, O.~Ogunlade,
  E.~Zhang, Accelerated high-resolution photoacoustic tomography via compressed
  sensing, Physics in Medicine \& Biology 61~(24) (2016) 8908.

\end{thebibliography}


\begin{thebibliography}{1}
\expandafter\ifx\csname url\endcsname\relax
  \def\url#1{\texttt{#1}}\fi
\expandafter\ifx\csname urlprefix\endcsname\relax\def\urlprefix{URL }\fi
\expandafter\ifx\csname href\endcsname\relax
  \def\href#1#2{#2} \def\path#1{#1}\fi

\bibitem{Negahban2012}
S.~Negahban, P.~Ravikumar, M.~J. Wainwright, B.~Yu, A unified framework for
  high-dimensional analysis of {M}-estimators with decomposable regularizers,
  Statistical Science 27~(4) (2012) 538--557.

\bibitem{Hunt2018}
X.~J. Hunt, P.~Reynaud-Bouret, V.~Rivoirard, L.~Sansonnet, R.~Willett, A
  data-dependent weighted {LASSO} under {P}oisson noise, IEEE Transactions on
  Information Theory 65~(3) (2018) 1589--1613.
\newblock \href {https://doi.org/10.1109/TIT.2018.2869578}
  {\path{doi:10.1109/TIT.2018.2869578}}.

\bibitem{Bernstein_wiki}
Bernstein inequalities (probability theory),
  \url{https://en.wikipedia.org/wiki/Bernstein_inequalities_(probability_theory)}.

\bibitem{BernoulliInequality_wiki}
Bernoulli's inequality,
  \url{https://en.wikipedia.org/wiki/Bernoulli's_inequality}.

\end{thebibliography}

\end{document}


\begin{frontmatter}


\title{Supplemental Material for `Performance Bounds for LASSO under Multiplicative LogNormal Noise: Applications to Pooled RT-PCR Testing'}
\author[label1]{Richeek Das}
\ead{richeek@cse.iitb.ac.in}
\affiliation[label1]{organization={Indian Institute of Technology Bombay},
            addressline={Department of Computer Science and Engineering}}
\affiliation[label3]{organization={Indian Institute of Technology Bombay},addressline={Department of Electrical Engineering}}


\author[label3]{Aaron Jerry Ninan}
\ead{aaronjerry12@gmail.com}

\author[label1]{Adithya Bhaskar}
\ead{adithyabhaskar@cse.iitb.ac.in}

\author[label1]{Ajit Rajwade}
\ead{ajitvr@cse.iitb.ac.in}





\end{frontmatter}


This supplemental material contains a brief introduction in \cref{sec:prelims}, statements of key propositions in \cref{sec:statements} and their proofs in \cref{apd:prop1}, and some additional experimental results in \cref{sec:addnl_results}.

\section{Preliminaries}\label{sec:prelims}
We reproduce expressions for the noise model from the main paper here below. In all of the following, $\boldsymbol{x^*} \in \mathbb{R}^p$ is the ground truth vector of viral loads, $\bm{y} \in \mathbb{R}^n$ is the measurement vector, $\bm{w} \in \mathbb{R}^m$ is the noise vector, and $\boldsymbol{A} \in \{0,1\}^{n \times p}$ is the pooling matrix, with entries drawn independently from $\textrm{Bernoulli}(q)$ where $q \in (0,1)$. We also have $n < p$. 
\begin{equation}\label{eqn:cs_prob_noise}
    y_j = \bm{A^j}\bm{x^*}(1 + q_a)^{w_j},
\end{equation}
where $q_a \in (0,1]$  is an amplification factor and $w_j \sim \mathcal{N}(0,\sigma^2)$ where $\sigma < 1$. Using a first order Taylor series expansion, we have:
\begin{equation}\label{eqn:noisemodel}
        y_j \approx \bm{A^j}\bm{x^*} + [\bm{A^j x^*}\ln(1+q_a)]w_j.
\end{equation}
As argued in the main paper, this approximation is good because $\sigma < 1$ and due to division by factors larger than 1 or multiplication by factors less than 1 in higher order terms. We also reproduce expressions for the surrogate sensing matrix $\bm{\tilde{A}}$ (for a Bernoulli $(q)$ sensing matrix $\bm{A}$) and surrogate measurement vector $\bm{\tilde{y}}$ from the main paper here below:
\begin{equation}
    \bm{\Tilde{A}} = \frac{\bm{A}}{\sqrt{nq(1-q)}} - \frac{q\bm{1}_{n\times 1}\bm{1}_{p\times 1}^T}{\sqrt{nq(1-q)}}.
\end{equation}
\begin{equation}
    \bm{\Tilde{y}} = \frac{1}{(n-1)\sqrt{nq(1-q)}}\left( n\bm{y} - \sum_{l=1}^{n}y_l\bm{1}_{n\times 1} \right).
\end{equation}

Given the surrogate sensing matrix $\bm{\Tilde{A}}$ and observations $\bm{\Tilde{y}}$, we redefine the \textsc{Lasso} and \textsc{WLasso} estimators via Equations \eqref{eqn:lassosurr} and \eqref{eqn:wlassosurr} respectively, as given below: 
\begin{equation}\label{eqn:lassosurr}
    \bm{\hat{x}^L} = \operatorname*{argmin}_{\bm{x} \in \mathbbm{R}^{p}} \|\bm{\Tilde{y}} - \bm{\Tilde{A}x}\|_{2}^{2} + \gamma \sum_{k=1}^{p}\beta|x_{k}|,
\end{equation}
where $\gamma > 2$ is a constant and $\beta > 0$ is a data-dependent scalar.
\begin{equation}\label{eqn:wlassosurr}
    \bm{\hat{x}^{WL}} = \operatorname*{argmin}_{\bm{x} \in \mathbbm{R}^{p}} \|\bm{\Tilde{y}} - \bm{\Tilde{A}x}\|_{2}^{2} + \gamma \sum_{k=1}^{p}\beta_k |x_{k}|,
\end{equation}
where $\gamma > 2$ is a constant and $\beta_k > 0$ is the $k$th data-dependent scalar. The exact expressions of $\beta$ and $\beta_k$ depend on the observed data and form the crux of our work.

We want the \emph{smallest} possible values of $\beta_k$ which would faithfully satisfy the Condition $\mathscr{C}1$ (defined below) with high probability. This serves as our primary objective in deriving the data-dependent weights for our noise model 
\begin{equation}\label{eqn:asswts}
\textrm{Condition } \mathscr{C}1 \textrm{ on } \left(\{\beta_k\}_{k}\right):  |(\bm{A}^T(\bm{y} - \bm{Ax^*}))_{k}| \leq \beta_k \quad \forall  k=1,...,p.
\end{equation}

\section{Statements of the Propositions}
\label{sec:statements}
\begin{proposition} 
{\normalfont\textbf{Weight for \textsc{Lasso}:}}
\label{prop1}

\noindent
Consider $\bm{A},\bm{\bar{R}} \in \mathbbm{R}^{n \times p}$ and define the following:
\begin{equation*}
    W := \max_{i,j} \left(\bm{\bar{R}}^T \bm{A} \right)_{ij} \textrm{ for } 1 \leq i \leq p, 1 \leq j \leq p, \quad \text{and} \quad \kappa := \sigma \ln (1+q_a)
\end{equation*}%
where $\bar{R}_{k,l} = \left(\frac{na_{l,k} - \sum_{l'=1}^{n}a_{l',k}}{n(n-1)q(1-q)}\right)^2$ for $k = 1, \cdots, p; l = 1,\cdots, n$.
\\
\noindent
Consider the following assumptions:
\begin{center}
\begin{itemize}
\addtolength{\itemindent}{1em}
    \item[\textit{A1.}] $nq \geq 12\max(q,1-q)\ln(p)$
    \item[\textit{A2.}] $p \geq 2$
    \item[\textit{A3.}] $\sigma < \frac{1}{\sqrt{2}\ln \left( 1 + q_a \right)} \left(\ln\left(\frac{1}{1 - \left( 1 - \frac{1}{p^3} \right)^{1/n}} \right)\right)^{-1/2}$
    \item[\textit{A4.}] $\bm{A}$ \textit{is a Bernoulli Sensing Matrix}.
\end{itemize}    
\end{center}
\noindent
If the aforementioned assumptions hold, then there exist positive constants $c$, $c'$ such that with probability larger than $1 - \frac{c'}{p}$, the choice%
\begin{equation*}
    \beta := \kappa \hat{\Lambda} \sqrt{6W\ln(p)} + c\left( \frac{3\ln(p)}{n} + \frac{9\max(q^2, (1-q)^2)}{n^2q(1-q)}\ln(p)^2 \right) \hat{\Lambda}
\end{equation*}%
satisfies Condition $\mathscr{C}1$ (defined earlier in \cref{eqn:asswts}), where $\hat{\Lambda}$ is an estimator of $\lVert \bm{x^*} \rVert_1$ given by
\begin{equation*}
    \hat{\Lambda} := \frac{\sum_{i=1}^{n} y_i + \sqrt{\sum_{i=1}^{n} y_i^2} \frac{\kappa\sqrt{6\ln(p)}}{1 - \kappa\sqrt{2g(3\ln(p))}}}{nq - \sqrt{6nq(1-q)\ln(p)} - \max(q, 1-q)\ln(p)},
\end{equation*}
where for any $\theta \in \mathbb{R}$, we define
\begin{equation*}
    g(\theta) := \ln\left( \frac{1}{1 - \left( 1 - e^{-\theta} \right)^{1/n}} \right).
\end{equation*}
Furthermore, $c=126$ works as long as $n \geq 20$.
\end{proposition}

\bigskip
\begin{proposition} {\normalfont\textbf{ Weights for \textsc{WLasso}}}\label{prop2}

\noindent
With the same notations and assumptions as in Proposition 1, there exist positive constants $c$, $c'$ such that with probability larger than $1 - \frac{c'}{p}$, the choice (depending on $k$)
\begin{equation*}
    \beta_k := \sqrt{\bm{\bar{R}_{k}}^T \bm{y_2}}\frac{\kappa\sqrt{6\ln(p)}}{1 - \kappa\sqrt{2g(3\ln(p))}} + c\left( \frac{3\ln(p)}{n} + \frac{9\max(q^2, (1-q)^2)}{n^2q(1-q)}\ln(p)^2 \right) \hat{\Lambda}
\end{equation*}
\noindent
satisfies Condition $\mathscr{C}1$ (see \cref{eqn:asswts}), where $\bm{y_2}:=\bm{y}\odot\bm{y}$. Furthermore, $c=126$ works as long as $n \geq 20$.
\end{proposition}

\section{Proof of Proposition 1 and Proposition 2}\label{apd:prop1}
To choose data-dependent weights, we need to satisfy the Condition $\mathscr{C}1$. From \cref{eqn:asswts}, we know that if each $\beta_k$ in the \textsc{WLasso} is close to $\left|\left(\bm{\Tilde{A}}^T\left(\bm{\Tilde{y}} - \bm{\Tilde{A}x^*}\right)\right)_{k}\right|$ then the risk bounds on both \textsc{Lasso} and \textsc{WLasso} approach the oracle least squares estimator \cite{Negahban2012}, \cite[Sec. I-B]{Hunt2018}. Hence, our target is to find the smallest possible $\beta_k$'s that satisfy the Condition $\mathscr{C}1$ in \cref{eqn:asswts}. To do this, we derive high-probability concentration inequalities to upper bound the LHS of \cref{eqn:asswts}. The rest of the proof here explains how we proceed with this.

\noindent
For all vectors $\bm{r} \in \mathbbm{R}^n$, we can deduce the following starting from \cref{eqn:noisemodel} as
follows:
\begin{eqnarray}\label{eq:noisemodelR}
\sum_j r_j y_j = \sum r_j \bm{A^j x^*} + \sum_j r_j \bm{A^j x^*} \ln (1+q_a) w_j \nonumber \\
\implies    \bm{r}^{T}\bm{y} \sim \mathcal{N}\left[ \bm{r}^T\bm{Ax^*}, \bm{r_2}^{T}(\bm{Ax^*})_{2} \kappa^2 \right],
\end{eqnarray}
 where $\kappa := \sigma\ln(1+q_a)$ and $\bm{r_2} := \bm{r}\odot\bm{r}$, that is the element-wise square of each element in the vector $\bm{r}$. The implication above follows because of the mutual independence of $\{w_j\}_{j=1}^n$.

\noindent
Now let us rewrite $\left(\bm{\Tilde{A}}^T\left(\bm{\Tilde{y}} - \bm{\Tilde{A}x^*}\right)\right)$ as:
\begin{align}\label{T1T2}
    \bm{\Tilde{A}^T}(\bm{\Tilde{y}} - \bm{\Tilde{A}x^*}) = {} & \bm{t_1} + \bm{t_2} \\
    = {} & \underbrace{\frac{1}{nq(1-q)}\left( \frac{n}{n-1}\bm{A}^T - \frac{\bm{A}^T\bm{1}_{n\times 1}\bm{1}_{n\times 1}^{T}}{n-1} \right)(\bm{y} - \bm{Ax^*})}_{\boldsymbol{t_1}} + \\
    & \frac{1}{nq(1-q)}\left( \frac{\bm{A}^T\bm{Ax^*} - \bm{A}^T\bm{1}_{n\times 1}\bm{1}_{n \times 1}^{T}\bm{Ax^*}}{n-1} \right. \nonumber \\
    & \underbrace{\left. + q\lVert \bm{x^*} \rVert_1 \bm{A}^T\bm{1}_{n\times 1} + q\bm{1}_{p\times 1}\bm{1}_{n \times 1}^{T}\bm{Ax^*} -q^2n\lVert \bm{x^*} \rVert_1 \bm{1}_{p \times 1} \vphantom{\frac{A^TAx^* - A^T\bm{1}_{n\times 1}\bm{1}_{n \times 1}^{T}Ax^*}{n-1}} \right)}_{\boldsymbol{t_2}}
\end{align}
which is easy to observe since $\sum_{l=1}^{n}\bm{y}_l\bm{1}_{n\times 1} = \bm{1}_{n \times 1}\bm{1}_{1 \times n} \bm{y}$ and $\bm{1}_{1 \times n}\bm{1}_{1 \times n}^T = n$.

\noindent
Then we have:
\begin{gather}
    (\bm{\Tilde{A}^T}(\bm{\Tilde{y}} - \bm{\Tilde{A}x^*}))_{k} = \bm{R_{k}}^{T}(\bm{y} - \bm{Ax^{*}}) + b_{k}(\bm{A}, \bm{x^{*}})
\end{gather}

\noindent
where $b_k (\bm{A},\bm{x^*})$ is defined below, and $\bm{R_{k}} \in \mathbbm{R}^n$ is the \textit{k}-th column of the matrix $\bm{R} \in \mathbbm{R}^{n \times p}$ as defined below:
\begin{gather}
\label{eq:Rb_k}
    \bm{R}^T := \frac{1}{nq(1-q)}\left( \frac{n}{n-1}\bm{A}^T - \frac{\bm{A}^T\bm{1}_{n\times 1}\bm{1}_{n\times 1}^{T}}{n-1} \right), \\
    b_k (\bm{A},\bm{x^*}) := t_{2,k}; t_{1k} = \bm{R_{k}}^{T}(\bm{y} - \bm{Ax^{*}}).
\end{gather}

Our aim is to upper bound $\|\bm{\Tilde{A}^T}(\bm{\Tilde{y}} - \bm{\Tilde{A}x^*})\|_{\infty}$, for which we need to upper bound $\lVert \bm{t_1} \rVert_{\infty} + \lVert \bm{t_2} \rVert_{\infty}$ in order to derive the constant weight in \cref{prop1}. To derive the non-constant weights in \cref{prop2} we set bounds on $| t_{1,k} | + | t_{2,k} |$. In the initial stages, we will see that all these bounds will depend on the input signal $\bm{x^*}$ or some quantities derived from it. Later on, we will present bounds and machinery to replace terms directly involving $\bm{x^*}$ with computable quantities.



Since we are using Bernoulli sensing matrices, we can reuse the calculations presented in \cite[Appendix C]{Hunt2018}. From \cite[Appendix C]{Hunt2018}, we state the probabilistic upper bound on $\|\bm{t_2}\|_{\infty}$ with some assumptions in the following \cref{T2interbound}:

\begin{lemma}\label{T2interbound}
Consider positive constants $c_1$ and $c_2$. Given $n \geq 2$ and $\theta > 1$, with probability larger than $1 - c_1 p^2 e^{-\theta}$, we have:
\begin{equation}
    \lVert \bm{t_2} \rVert_{\infty} \leq c_2 \left( \frac{\theta}{n} + \frac{\max(q^2, (1-q)^2)}{n^2q(1-q)}\theta^2 \right) \lVert \bm{x^*} \rVert_1, 
\end{equation}
where $c_2=126$ works when $n \geq 20$.
\end{lemma}

\medskip
The bound on $\lVert \bm{t_2} \rVert_{\infty}$ in \cref{T2interbound}, requires the knowledge of $\lVert \bm{x^*} \rVert_1$. Since it is generally impossible to have knowledge of the $l_1$ norm of the input signal in a reconstruction problem, we derive an estimator for it later in the proof. 

\noindent At first, we develop the required intermediate bounds and estimators needed to set concentration bounds on $|t_{1,k}|$. From Equation \eqref{T1T2} we observe:
    \begin{equation}\label{eq:T1status}
        t_{1,k} = \bm{R_k}^T(\bm{y}-\bm{Ax^*}), 
    \end{equation}
    where $\bm{R}$ is as defined in Equation~\ref{eq:Rb_k} and $\bm{R_k}$ is its $k$th column. That is, for all $l=1,\ldots, n$, and $k=1,\cdots,p$ we have:
    \begin{equation}
        R_{k,l} = \frac{na_{l,k} - \sum_{l'=1}^{n}a_{l',k}}{n(n-1)q(1-q)}.
    \end{equation}

\noindent To bound $\|\bm{t_1}\|_{\infty}$ and $|t_{1,k}|$, we will need some concentration inequalities on $\|\bm{R}^T(\bm{y}-\bm{Ax^*})\|_{\infty}$. With this goal in mind, we now develop the required concentration inequalities for our case of multiplicative lognormal noise.
    
\subsection{Concentration Bounds for Multiplicative Lognormal Noise}

\noindent With the same terminology presented in \cref{eq:noisemodelR}, we follow the lines of the proof of Bernstein's inequality. Let $z := \bm{r}^T\bm{y}$ and $\Tilde{z} := \bm{r}^T\bm{Ax^*}$, where $\bm{r} \in \mathbbm{R}^{n}$. Therefore for all $\lambda \in \mathbbm{R}: \lambda > 0$ we have:
\begin{gather*}
        \mathbbm{E}\left( e^{\lambda (z-\Tilde{z})} | \bm{A} \right) =  \mathbbm{E}\left( e^{\lambda \bm{r}^T [\bm{y} - \bm{Ax^*}]} | \bm{A} \right) 
        = \mathbbm{E}\left( e^{\lambda \bm{r}^T\bm{y}} | \bm{A} \right) e^{-\lambda \bm{r}^T\bm{Ax^*}}.
\end{gather*}
Using the standard expression for the moment generating function (MGF) of $\bm{r}^T\bm{y}$, which is Gaussian distributed as per \cref{eq:noisemodelR}, we have:
\begin{gather*}
\mathbbm{E}\left( e^{\lambda (z-\Tilde{z})} | \bm{A} \right)        = \exp\left(\bm{r}^T\bm{Ax^*}\lambda + \frac{1}{2}\bm{r_{2}}^T(\bm{Ax^*})_2\kappa^2\lambda^2\right) \exp{(-\lambda \bm{r}^T\bm{Ax^*})} \\
        = \exp\left(\frac{1}{2}\bm{r_2}^T(\bm{Ax^*})_2\kappa^2\lambda^2\right),
\end{gather*}
where $\bm{r_2} := \bm{r} \odot \bm{r}$ as defined before and $\lambda$ is the MGF parameter. 
Therefore,
    \begin{equation}
        \mathbbm{E}\left( e^{\lambda (z-\Tilde{z})} | \bm{A} \right) = \exp\left( \frac{\kappa^2\lambda^2 v}{2} \right), \text{ where } v := \bm{r_2}^T(\bm{Ax^*})_2.
    \end{equation}
Using Markov's Inequality, for all $u \in \mathbbm{R} : u>0$
    \begin{equation}
        \mathbbm{P}(z-\Tilde{z} \geq u) \leq \exp\left( \frac{\kappa^2\lambda^2v}{2} - \lambda u \right).
    \end{equation}
    We need to find the minimal upper bound, so we optimize on $\lambda$, thus yielding $\lambda = \frac{u}{\kappa^2v}$.
    This gives us
    \begin{equation}
      \mathbbm{P}\left(z-\Tilde{z} \geq u\right) \leq \exp\left( -\frac{u^2}{2\kappa^2v} \right).
    \end{equation}
    Denoting the upper bound on the probability by $e^{-\theta}$, we obtain $u = \kappa \sqrt{2\theta v}$.

\noindent
Therefore we have our first bound
    \begin{equation}\label{1stbound}
        \mathbbm{P}\left(\bm{r}^T\bm{y} - \bm{r}^T\bm{Ax^*} \geq \kappa \sqrt{2\theta v}\right) \leq e^{-\theta}. 
    \end{equation}

    
\noindent
Replacing $\bm{r}$ by $-\bm{r}$ and using the symmetry of a Gaussian distribution, we obtain the following via a union bound:
\begin{equation}\label{2ndbound}
    \mathbbm{P}\left(\left| \bm{r}^T\bm{y} - \bm{r}^T\bm{Ax^*} \right| \geq \kappa \sqrt{2\theta v}\right) \leq 2e^{-\theta}. 
\end{equation}

To proceed with the proof, we need to find bounds for $|\bm{R_l}^T(\bm{y}-\bm{Ax^*})|$ such that they are independent of the input signal $\bm{x^*}$. In our case, the above bounds in \cref{1stbound} and \cref{2ndbound} still depend on $v := \bm{r}_2^T(\bm{Ax^*})_2$. We still need to find bounds for $v$.
     
\subsection{Bounding $v = \bm{r}_2^T(\bm{Ax^*})_2$}

From \cref{eqn:noisemodel}, we obtain the following result when we set $\eta_i := w_i/\sigma$ where $w_i \sim \mathcal{N}(0,\sigma^2)$: 
\begin{equation}
    (\bm{Ax^*})_{i}^{2} = \dfrac{\bm{y}_i^2}{(1 + \eta_{i}\kappa)^2}.
\end{equation}
By definition of $\eta_i$, we have $\eta_i \sim \mathcal{N}(0,1)$. Hence, we can use a standard Chernoff bound to obtain an upper concentration bound on $\eta_i$ with any $\varphi \in \mathbbm{R}$:
\begin{equation}
    \mathbbm{P}\left( \eta_{i} \leq -\frac{\varphi}{\kappa} \right) \leq \exp\left(-\frac{\varphi^{2}}{2\kappa^{2}}\right).
\end{equation}

Therefore under the assumption $\varphi < 1$, with probability $\geq 1 - \exp\left(-\frac{\varphi^{2}}{2\kappa^{2}}\right)$, we have
\begin{equation}
    \frac{\bm{y}_i^2}{(1 + \eta_{i}\kappa)^2} < \frac{\bm{y}_i^2}{(1 - \varphi^{2})}.
\end{equation}
Hence, we have:
\begin{equation}
 \mathbbm{P}\left( \bm{r}_{i}^2(\bm{Ax^*})_i^2 \geq \frac{\bm{r}_{i}^{2}\bm{y}_{i}^2}{(1 - \varphi)^2}\right) \leq \exp\left(-\frac{\varphi^{2}}{2\kappa^{2}}\right) 
\end{equation}
\begin{equation}
    \implies \mathbbm{P}\left( \bm{r}_{i}^2(\bm{Ax^*})_i^2 < \frac{\bm{r}_{i}^{2}\bm{y}_{i}^2}{(1 - \varphi)^2}\right) > 1 - \exp\left(-\frac{\varphi^{2}}{2\kappa^{2}}\right) 
\end{equation}
The event $\bm{r_2}^{T} (\bm{Ax^*})_{2} < \frac{\bm{r_2}^{T}\bm{y_2}}{(1 - \varphi)^2}$ will hold true when each of the events $\mathscr{E}_i := \bm{r}_{i}^2(\bm{Ax^*})_i^2 < \frac{\bm{r}_{i}^{2}\bm{y}_{i}^2}{(1 - \varphi)^2}$ (across $i$) occurs. Given the mutual independence of $w_i$'s, we know that the events $\{\mathscr{E}_i\}_{i=1}^n$ are independent. Hence, we obtain:
\begin{equation}
     \mathbbm{P} \left( \bm{r_2}^{T} (\bm{Ax^*})_{2} < \frac{\bm{r_2}^{T}\bm{y_2}}{(1 - \varphi)^2} \right) > \left(1 - \exp{\left(-\frac{\varphi^2}{2\kappa^2}\right)}\right)^n.
\end{equation}
\begin{equation}
    \implies \mathbbm{P} \left( \bm{r_2}^{T} (\bm{Ax^*})_{2} \geq \frac{\bm{r_2}^{T}\bm{y_2}}{(1 - \varphi)^2} \right) \leq 1 - \left(1 - \exp{\left(-\frac{\varphi^2}{2\kappa^2}\right)}\right)^n.
\end{equation}

\medskip
\noindent
Setting $\exp(-\theta) = 1 - \left(1 - \exp{\left(\frac{-\varphi^2}{2\kappa^2}\right)}\right)^n$, we obtain:      
        \begin{equation*}
        \label{kgthetaassump}
            \varphi = \kappa \sqrt{2\ln\left( \frac{1}{1 - \left( 1 - e^{-\theta} \right)^{1/n}} \right)}.
        \end{equation*}
        
        Since $\varphi < 1$, the earlier equation can be restated more concisely in the form of a new function $g(\theta)$:
        \begin{equation}\label{sigmapassumption}
            g(\theta) = \ln\left( \frac{1}{1 - \left( 1 - e^{-\theta} \right)^{1/n}} \right) \quad \text{and} \quad \kappa\sqrt{2g(\theta)} < 1.
        \end{equation}
        \textbf{Note: } Expanding \cref{sigmapassumption} leads to our Assumption \textit{A3} in \cref{prop1} and \cref{prop2}. 

\medskip
Under the above assumption, we have our high-probability upper bound on $v = \bm{r_2}^{T} (\bm{Ax^*})_{2}$ as:
        \begin{equation}
        \label{eq:r2ax2}
            \mathbbm{P}\left( \bm{r_2}^{T} (\bm{Ax^*})_{2} \geq \frac{\bm{r_2}^{T}\bm{y_2}}{(1 - \kappa\sqrt{2g(\theta)})^2} \right) \leq e^{-\theta}
        \end{equation}

\noindent
Combining \cref{1stbound}, \cref{2ndbound} and \cref{eq:r2ax2} we obtain the required concentration inequalities on $|\bm{R_l^T}(\bm{y}-\bm{Ax^*})|$. For instance, \cref{1stbound} and \cref{eq:r2ax2} can be easily combined through a union bound on the following two events: 
\begin{gather}
\mathscr{X}_1 := \bm{r}^T\bm{y} - \bm{r}^T\bm{Ax^*} \geq \kappa \sqrt{2\theta \bm{r_2}^{T} (\bm{Ax^*})_{2}} \\
\mathscr{X}_2 := \bm{r_2}^{T} (\bm{Ax^*})_{2} \geq \frac{\bm{r_2}^{T}\bm{y_2}}{(1 - \kappa\sqrt{2g(\theta)})^2},   
\end{gather}
since our target event $\mathscr{X}_3 := \left(\bm{r}^T\bm{y} - \bm{r}^T\bm{Ax^*} \geq \sqrt{\bm{r_2}^{T}\bm{y_2}} \frac{\kappa \sqrt{2\theta}}{1 - \kappa\sqrt{2g(\theta)}} \right)$ might hold true even if one of events $\mathscr{X}_1$ or $\mathscr{X}_2$ does not occur. We state this formally in the following \cref{lemma:ryax}.
\begin{lemma}\label{lemma:ryax}
Let $\bm{R} \in \mathbbm{R}^{n \times p}$ be any real-valued matrix and $\bm{\bar{R}} = \bm{R} \odot \bm{R}$ be a matrix containing the element-wise squares of $\bm{R}$, then the following high-probability upper bounds hold
        \begin{equation}
        \label{eq:ryax}
            \mathbbm{P}\left( \bm{R_l}^{T} (\bm{y} - \bm{Ax^*}) \geq \sqrt{\bm{\bar{R}_l}^{T}\bm{y_2}} \frac{\kappa \sqrt{2\theta}}{1 - \kappa\sqrt{2g(\theta)}} \right) \leq 2e^{-\theta}
        \end{equation}
        \begin{equation}
        \label{eq:ryaxmod}
            \mathbbm{P}\left( \left| \bm{R_l}^{T} (\bm{y} - \bm{Ax^*}) \right| \geq \sqrt{\bm{\bar{R}_l}^{T}\bm{y_2}} \frac{\kappa \sqrt{2\theta}}{1 - \kappa\sqrt{2g(\theta)}} \right) \leq 3e^{-\theta}
        \end{equation}
under the assumption $\kappa\sqrt{2g(\theta)} < 1$ where,
    \begin{equation*}
        g(\theta) = \ln\left( \frac{1}{1 - \left( 1 - e^{-\theta} \right)^{1/n}} \right) \quad \text{and} \quad \kappa := \ln(1+q_a)\sigma
    \end{equation*}
\end{lemma}
        

\subsection{Bounding $\left| t_{1,k} \right|$ and $\lVert \bm{t_1} \rVert_{\infty}$}

With \cref{lemma:ryax} in place, we can set the required bounds on $\left| t_{1,k} \right|$ and $\lVert \bm{t_1} \rVert_{\infty}$ for the non-constant (weighted) and constant (plain) weights respectively. From \cref{eq:T1status} we know:
\begin{gather*}\label{eqn:t1kbound}
    t_{1,k} = \bm{R_k}^T(\bm{y} - \bm{Ax^*}) \quad and \quad
    R_{k,l} = \frac{na_{l,k} - \sum_{l'=1}^{n}a_{l',k}}{n(n-1)q(1-q)}
\end{gather*}


\noindent
\textbf{Note:} $\boldsymbol{R} \in \mathbbm{R}^{n\times p}$ and $\boldsymbol{R_k} \in \mathbbm{R}^{n}$, $k = 1, \cdots, p$ and $l = 1, \cdots, n$.

Define, $\bar{R}_{k, l} = R_{k, l}^2$, i.e. $\bm{\bar{R}} = \bm{R}\odot\bm{R}$. By \cref{eq:ryaxmod} of \cref{lemma:ryax}, we can directly obtain a bound for $|t_{1,k}|$.  We have the following result with probability $\geq 1 - 3e^{-\theta}$:
    \begin{equation}
    \label{t1kboundmain}
        \left| t_{1,k} \right| \leq \sqrt{\bm{\bar{R}_k}^T\bm{y_2}} \frac{\kappa\sqrt{2\theta}}{1 - \kappa\sqrt{2g(\theta)}}.
    \end{equation}
For constant weights we need to find $\lVert \boldsymbol{t_1} \rVert_{\infty}$. Before proceeding in this direction, we will first establish the following auxiliary lemma:
\begin{lemma}\label{lemma:randombound}
For any vectors $\bm{s} \in \mathbbm{R}^{+p}$, $\bm{x} \in \mathbb{R}^p$ and $a_{l,m} \in \{0, 1\}$, we have
\begin{equation}
    \sum_{l=1}^n s_l \left(\sum_{m=1}^p a_{l,m}x_m\right)^2 \leq \max_{m \in \{1,2,\cdots,p\}} \left(\sum_{l=1}^n s_l a_{l,m}\right) \|\bm{x}\|_1^2.
    \label{eq:Lemma1} \blacksquare
\end{equation}
\noindent
\textit{\textbf{Proof.}} We can assume without loss of generality that $x_i \geq 0$ for all $i \in \{1,2,\cdots,p\}$, since if this were not true, we could always find an index $i$ with $x_i < 0$, whereupon the transformation $x_i \longrightarrow -x_i$ would increase the LHS but keep the RHS the same. A series of such transformations can always lead us to an inequality where all elements of $x_i$ are non-negative, which is stronger than the original result (and is proved here).
\noindent
Let $k$ be the index defined as
\begin{equation*}
    k \equiv \argmax_m \sum_{l=1}^n s_l a_{l,m} = \argmax_m \sum_{l:a_{l,m}=1} s_l
\end{equation*}
In other words, the index $k$ maximizes the sum of $s_l$'s chosen by the position of $1$'s in the column $k$. Then, note that for any two indices $u,v$, we have
\begin{equation*}
    \sum_{l=1}^n s_l a_{l,k} \geq \sum_{u: a_{l,u}=1} s_l \geq \sum_{u,v:a_{l,u}=1\atop a_{l,v}=1} s_l = \sum_{l=1}^n s_l a_{l,u}a_{l,v}
\end{equation*}
Following from this and the definition of $k$ we have,
\begin{align*}
\max_{m \in \{1,2,\cdots,p\}} \left(\sum_{l=1}^n s_l a_{l,m}\right) \|\bm{x}\|_1^2 &= \left(\sum_{l=1}^n s_l a_{l,k}\right) \cdot \left(\sum_{u=1}^p x_u\right) \cdot \left(\sum_{v=1}^p x_v\right) \\
&= \sum_{l\in \{1,2,\cdots,n\}\atop u,v \in \{1,2,\cdots,p\}} s_la_{l,k}x_ux_v \\
&\geq \sum_{l\in \{1,2,\cdots,n\}\atop u,v \in \{1,2,\cdots,p\}} s_la_{l,u}x_ua_{l,v}x_v \\
&= \sum_{l=1}^n s_l \left(\sum_{m=1}^p a_{l,m}x_m\right)^2
\end{align*}
as claimed. $\blacksquare$
\end{lemma}

\noindent
From Equation~\ref{2ndbound}, we have:
\begin{equation}
\label{t1kmod}
    \mathbbm{P} \left( \left| t_{1,k} \right| \geq \kappa \sqrt{2\theta \bm{\bar{R}_{k}}^{T}(\bm{Ax^*})_2} \right) \leq 2e^{-\theta}.
\end{equation}
The RHS involves a signal-dependent term $\bm{\bar{R}_k}^T(\bm{Ax^*})_{2}$. We want to upper bound this quantity, independent of $k$, to obtain a concentration bound on $\lVert \bm{t_1} \rVert_{\infty}$:
\begin{align}
\label{vkax2}
    {\bm{\bar{R}_k}^T(\bm{Ax^*})_{2}} &= \sum_{l=1}^{n} r_{k,l}^2 \left( \sum_{m=1}^{p} a_{l, m}x^*_m \right)^2 \nonumber \\
    &\leq \max_{k=1,\cdots,p\atop m=1,\cdots,p} \left(\sum_{l=1}^{n} r_{k,l}^2 a_{l,m}\right) \lVert \bm{x^*} \rVert_{1}^{2} \quad \textit{using \cref{lemma:randombound} with } s_k := r^2_{k,l} \nonumber \\
    &= \lVert \bm{x^*} \rVert_1^2 W,
\end{align}
where $W := \max_{k=1,\cdots,p \atop m=1,\cdots,p}\left(\sum_{l=1}^{n} r_{k,l}^2 a_{l,m}\right) = \max_{k=1,\cdots,p \atop m=1,\cdots,p}\left(\left( \bm{\bar{R}}^T\bm{A} \right)_{k,m}\right)$. 
    
\medskip
From \cref{t1kmod} and \cref{vkax2} we have, with probability $\geq 1 - 2pe^{-\theta}$,
\begin{equation}\label{eqn:t1Wx1}
    \lVert \bm{t_1} \rVert_{\infty} \leq \kappa \sqrt{2\theta\lVert \bm{x^*} \rVert_{1}^{2}W}.
\end{equation}
The probability $1 - 2pe^{-\theta}$ is obtained from a union bound. However, the quantity $\|\bm{x^*}\|_1$ is still unknown. We need to derive an estimator for $\|\bm{x^*}\|_1$ to bound $\lVert \bm{t_1} \rVert_{\infty}$ in computable terms.

\subsection{Estimating $\lVert \bm{x^*} \rVert_1$}


As mentioned in the discussions of both \cref{eqn:t1Wx1} and \cref{T2interbound}, to obtain an input signal independent upper bound on $\lVert \boldsymbol{t_1} \rVert_{\infty}$ and $\lVert \boldsymbol{t_2} \rVert_{\infty}$, we need to derive an estimator which presents theoretical guarantees for the upper bound of $\lVert \boldsymbol{x^*} \rVert_1$. 
Instead of a direct upper bound, it is much easier to set a bound on $\sum_{l,k} a_{l,k}x_{k}^{*}$ (recall that all elements of $\bm{x^*}$ are non-negative) and invoke a second inequality to obtain the concentration bound on $\lVert \bm{x^*} \rVert_1$.
With this in mind, since $\bm{A}$ is a Bernoulli Sensing Matrix with a parameter \textit{q}, we can invoke Bernstein's inequality. To make the paper self-contained, we state Bernstein's inequality from \cite{Bernstein_wiki} for our case of Bernoulli random variables here below  as \cref{lemma:bernsteinbernoulli}.
\noindent
\begin{lemma}
\label{lemma:bernsteinbernoulli}
\textbf{(Bernstein Inequality) } Let $X_1,\cdots,X_n$ be independent zero-mean random variables. Suppose that $|X_i|\leq M$ for all $i$, then for all positive $t$, we have 
\begin{equation}    \label{basicbernstein}\mathbbm{P}\left(\sum_{i=1}^nX_i \geq t\right) = \mathbbm{P}\left(\sum_{i=1}^nX_i \leq -t\right) \leq \exp\left(-\frac{\frac{1}{2}t^2}{\sum_{i=1}^n\mathbbm{E}[X_i^2]  + \frac{1}{3}Mt }\right).
\end{equation}
\noindent
We adapt this to our case of Bernoulli random variables. Let $S_n := \sum_{i=1}^nX_i$ where $X_i \sim \text{Ber}(q)$. Then,
\begin{equation}
\label{bernoullibernstein}
    \mathbbm{P}\left( S_n \leq nq - \frac{\theta}{3}\max(q,1-q) - \sqrt{2\theta nq(1-q)} \right) \leq e^{-\theta} .\blacksquare
\end{equation}

\noindent
\textit{\textbf{Proof:}} To apply Bernstein inequality, we first center $X_i \leftarrow X_i - \mathbbm{E}[X_i]$. Then, we have
\begin{equation*}
    \left| X_i - \mathbbm{E}[X_i] \right| = \left| 
X_i - q \right| \leq \max(q,1-q) \quad \text{and} \quad
\sum_{i=1}^n\mathbbm{E}[X_i^2] = nq(1-q).
\end{equation*}
Invoking \cref{basicbernstein} we can state,
\begin{equation}
\label{intermediatebernoulibern}
    \mathbbm{P}\left(S_n - nq \leq -t\right) \leq \exp\left( -\frac{t^2/2}{\frac{t}{3}\max(q,1-q) + nq(1-q)} \right).
\end{equation}
Let us define $\theta := \dfrac{t^2/2}{\frac{t}{3}\max(q,1-q) + nq(1-q)}$. Rewriting,
\begin{align*}
    t &= \frac{\theta}{3}\max(q,1-q) + \sqrt{\frac{\theta^2\left(\max(q,1-q)\right)^2}{9} + 2\theta nq(1-q)} \\
    &\geq \frac{\theta}{3}\max(q,1-q) + \sqrt{2\theta nq(1-q)}.
\end{align*}
Hence, this inequality combined with \cref{intermediatebernoulibern} gives us the required statement in \cref{bernoullibernstein}. $\blacksquare$
\end{lemma}

Since $a_{l,k}$ are i.i.d. Bernoulli random variables, such that $a_{l,k} \sim \text{Ber}(q)$, we can utilise \cref{lemma:bernsteinbernoulli} with $X_l = a_{l,k}$ where $l \in \{1,\cdots,n\}$. Therefore for any given $k=1,\cdots,p$, we have with probability larger than $1-e^{-\theta}$:
\begin{equation}
\label{x1xarel}
    \sum_{l=1}^n (a_{l,k} - q) \geq -C_{n,\theta} \quad where \quad C_{n,\theta} := \sqrt{2nq(1-q)\theta} + \max(q, 1-q)\frac{\theta}{3}.
\end{equation}%
Hence on the same event for any given $k$ we have:
\begin{equation*}
    nq - C_{n,\theta} \leq \sum_{l=1}^{n} a_{l,k} = \bm{1}_{1\times n}\bm{A}_k \implies \left(nq - C_{n,\theta}\right)x^*_k \leq \bm{1}_{1\times n}\bm{A}_k x^*_k.
\end{equation*}
Summing over $k=1,\cdots,p$ we have the following expression:
\begin{equation}\label{eqn:l1sumrel}
    (nq - C_{n,\theta})\lVert \bm{x^*} \rVert_1 \leq \bar{x}_{a} \quad \text{where} \quad \bar{x}_a = \sum_{l,k} a_{l,k}x_{k}^{*} = \bm{1}_{1 \times n} \bm{Ax^*},
\end{equation}%
which holds with probability $\geq (1-e^{-\theta})^p$, since the columns $\bm{A}_k$ are mutually independent. Now since $p \geq 1$ and $e^{-\theta} < 1$, we can state $(1-e^{-\theta})^p \geq 1-pe^{-\theta}$ via Bernoulli's inequality \cite{BernoulliInequality_wiki}. For simplicity we consider \cref{eqn:l1sumrel} to be true with probability $\geq 1-pe^{-\theta}$.

The first assumption on the range of $(n,q)$ is to assume that $nq > C_{n,\theta}$, which is implied by%
\begin{equation}
\label{qassumption1}
    nq \geq 4\max(q,1-q)\theta.
\end{equation}
from a similar argument as presented in Comment 3 of the discussion following \cref{prop1} and \cref{prop2}.

Now, to obtain an upper bound on $\lVert \bm{x^*} \rVert_1$, we need to find some concentration bound on $\bar{x}_a$. We bound $\boldsymbol{\bm{1}_{n \times 1} \bm{Ax^*}}$ using \cref{lemma:ryax}. Setting $\bm{R} = -\bm{1}_{n\times 1}$ in \cref{eq:ryax}, we have with probability $\geq 1 - 2e^{-\theta}$:
\begin{equation}\label{1ax}
    \underbrace{\bm{1}_{1 \times n} \bm{Ax^*}}_{\bar{x}_a} < \sum_{i=1}^{n} y_i + \sqrt{\sum_{i=1}^{n} y_i^2} \frac{\kappa\sqrt{2\theta}}{1 - \kappa\sqrt{2g(\theta)}}.
\end{equation}

\noindent
In this case with probability $\geq 1 - (p+2)e^{-\theta}$

\begin{equation}\label{x1bound}
    \lVert \bm{x^*} \rVert_1 \leq \frac{1}{nq - C_{n,\theta}} \left( \sum_{i=1}^{n} y_i + \sqrt{\sum_{i=1}^{n} y_i^2} \frac{\kappa\sqrt{2\theta}}{1 - \kappa\sqrt{2g(\theta)}} \right) := \hat{\Lambda}_{\theta}.
\end{equation}
Note that $\hat{\Lambda}$ in \cref{prop1} is obtained from $\hat{\Lambda}_{\theta}$ defined above by setting $\theta := 3 \ln p$. We show a numerical illustration of the expression for $\hat{\Lambda}$ on over 1000 simulations in the section on numerical experiments in the main paper (Figure 4). Now that we have the required estimator $\hat{\Lambda}$ for $\|\bm{x^*}\|_1$, we can use it in \cref{eqn:t1Wx1} to obtain the required bound on $\|\bm{t_1}\|_{\infty}$. With probability $\geq 1 - (3p+2)e^{-\theta}$, we have:
    \begin{equation}
    \label{t1infty}
        \lVert \bm{t_1} \rVert_{\infty} \leq \kappa \hat{\Lambda}_{\theta} \sqrt{2\theta W}
    \end{equation}
Further combining \cref{T2interbound} and \cref{x1bound} along with their underlying assumptions, we state the required upper bound on $\lVert \bm{t_2} \rVert_{\infty}$ in the following \cref{lemma:T2}.

\begin{lemma}\label{lemma:T2}
Let there exist some positive constants $c_2$ and $c_3$ such that as soon as $n \geq 2$ and $\theta > 1$, with probability larger than $1 - c_3 p^2 e^{-\theta}$
\begin{equation}\label{t2bound}
    \lVert \bm{t_2} \rVert_{\infty} \leq c_2 \left( \frac{\theta}{n} + \frac{\max(q^2, (1-q)^2)}{n^2q(1-q)}\theta^2 \right) \hat{\Lambda}_{\theta},
\end{equation}
where
\begin{equation*}
    \hat{\Lambda}_{\theta} := \frac{1}{nq - C_{n,\theta}} \left( \sum_{i=1}^{n} y_i + \sqrt{\sum_{i=1}^{n} y_i^2} \frac{\kappa\sqrt{2\theta}}{1 - \kappa\sqrt{2g(\theta)}} \right),
\end{equation*}
and
\begin{equation*}
    C_{n,\theta} = \sqrt{2nq(1-q)\theta} + \max(q, 1-q)\frac{\theta}{3},
\end{equation*}
under the assumption $nq \geq 4\max(q,1-q)\theta$. 
\end{lemma}

\medskip
\textbf{Note:} \cref{t1infty} combined with \cref{lemma:T2} for $\|\bm{t_2}\|_{\infty}$ produces the constant weight stated in \cref{prop1}. On the other hand, \cref{t1kboundmain} combined with \cref{lemma:T2} leads to \cref{prop2}. In both cases, we set $\theta := 3\ln{p}$ which gives us bounds having high probability $\geq 1 - \frac{c'}{p}$.

\begin{figure}[t!]
    \centering 
\begin{subfigure}{0.33\textwidth}
  \includegraphics[width=1.1\linewidth]{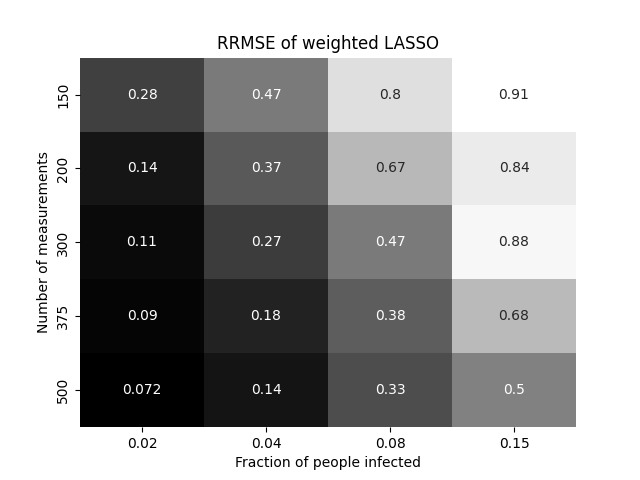}
\end{subfigure}\hfil 
\begin{subfigure}{0.33\textwidth}
  \includegraphics[width=1.1\linewidth]{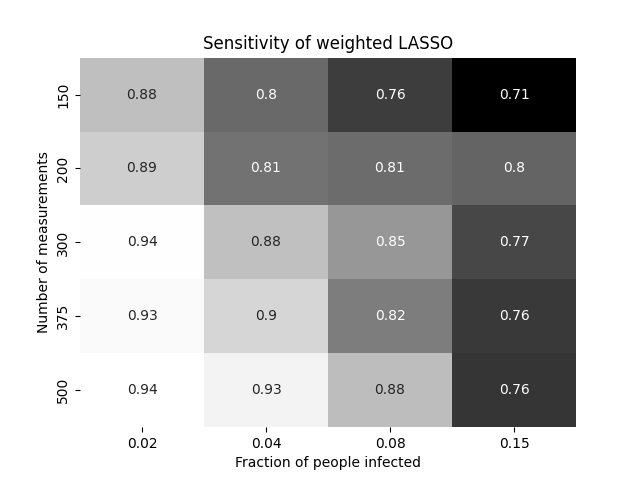}
\end{subfigure}\hfil 
\begin{subfigure}{0.33\textwidth}
  \includegraphics[width=1.1\linewidth]{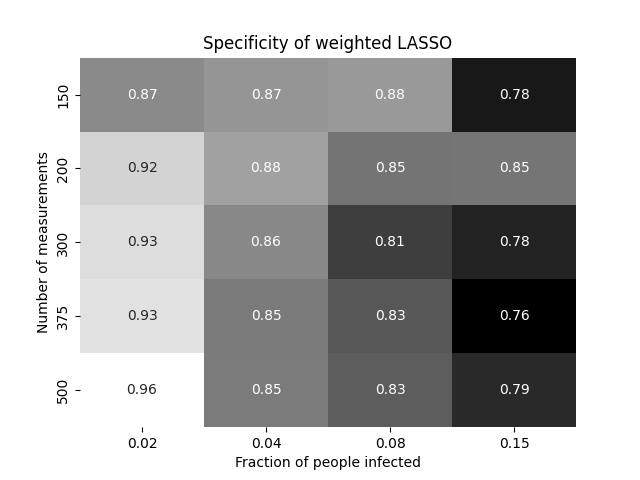}
\end{subfigure}

\begin{subfigure}{0.33\textwidth}
  \includegraphics[width=1.1\linewidth]{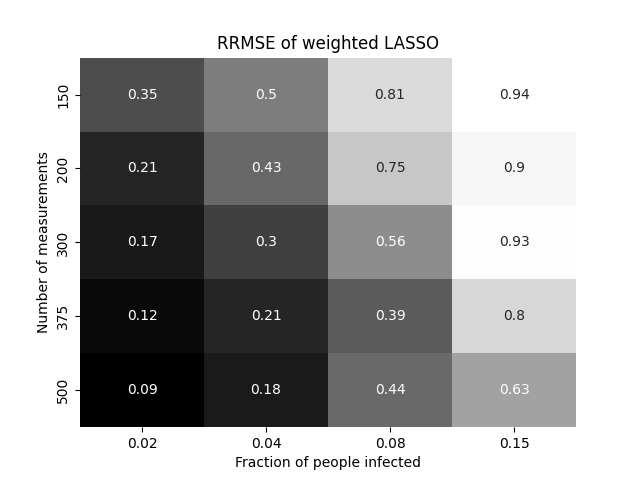}
\end{subfigure}\hfil
\begin{subfigure}{0.33\textwidth}
  \includegraphics[width=1.1\linewidth]{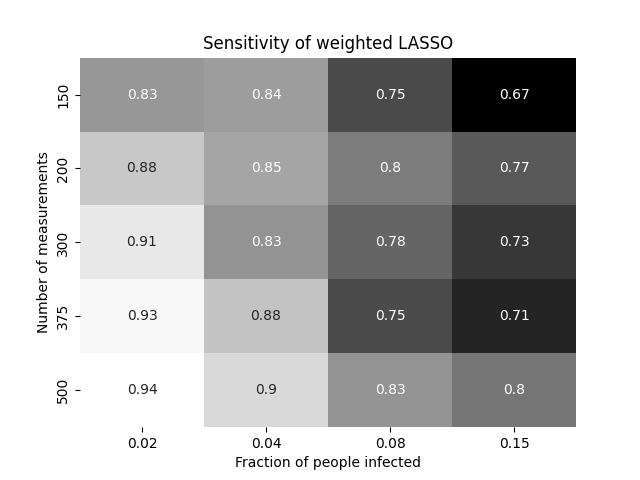}
\end{subfigure}\hfil
\begin{subfigure}{0.33\textwidth}
  \includegraphics[width=1.1\linewidth]{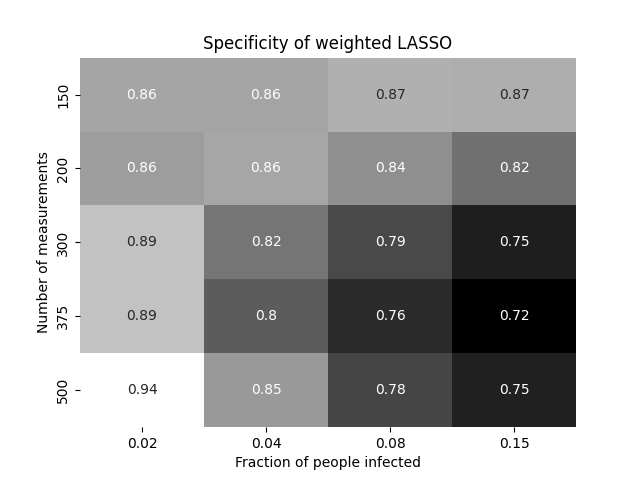}
\end{subfigure}

\begin{subfigure}{0.33\textwidth}
  \includegraphics[width=1.1\linewidth]{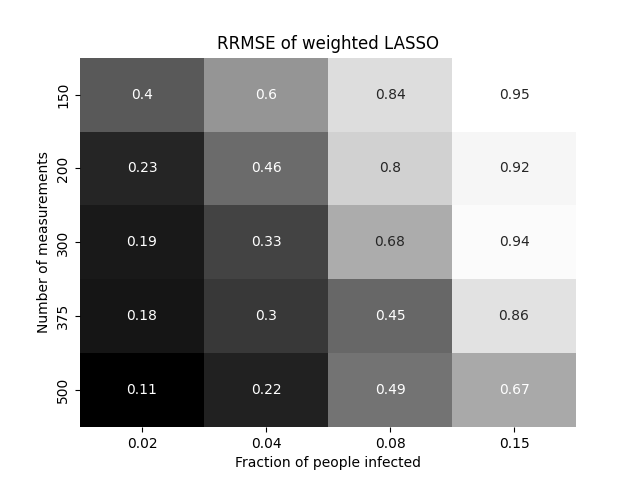}
\end{subfigure}\hfil
\begin{subfigure}{0.33\textwidth}
  \includegraphics[width=1.1\linewidth]{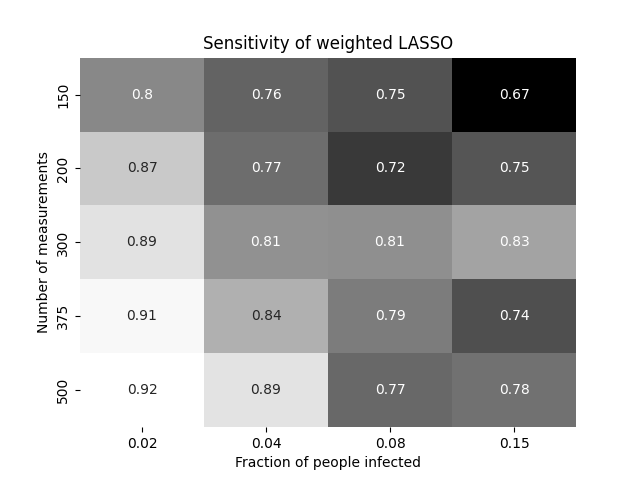}
\end{subfigure}\hfil
\begin{subfigure}{0.33\textwidth}
  \includegraphics[width=1.1\linewidth]{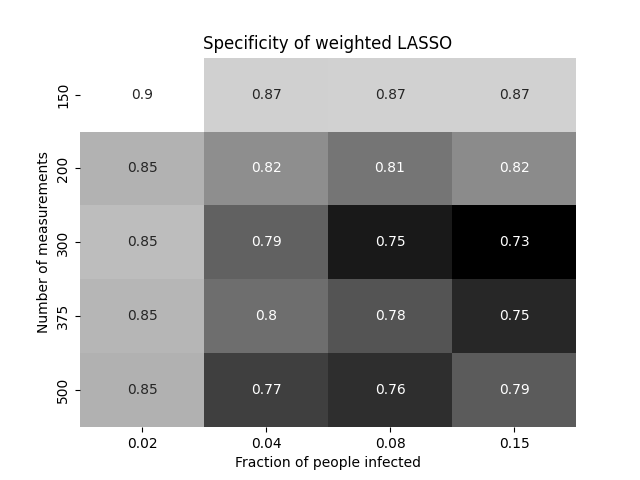}
\end{subfigure}
\caption{RRMSE, Sensitivity and Specificity values for $\sigma \in \{0.1,0.15,0.2\}$ (top to bottom), versus the number of measurements $n$ and $f_s$ (fraction of people infected out of a total of $p$). The results are for the \textsc{WLasso} algorithm. Lower is better for RRMSE, and higher is bette for the rest.}
\label{fig:results_diff_sigma}
\end{figure}

\section{Additional Numerical Results}
\label{sec:addnl_results}
We present some additional numerical results for the same experimental setup and parameters as in the main paper, except that we select $\sigma \in \{0.1,0.15,0.2\}$ over and above the experiment with $\sigma = 0.05$ in the main paper. Here $\sigma$ is the standard deviation of the Gaussian noise on noting cycle time values in RT-PCR. These results (RRMSE, sensitivity, specificity) are shown in Fig.~\ref{fig:results_diff_sigma}. The general trend of performance improving with increase in $n$ or decrease in $\ell_0$ norm ($ = f_s p$) are visible in these plots as well. Note that values of $\sigma$ higher than 0.2 are neither physically realistic nor would they satisfy Assumption A3 of Proposition 1 of the main paper. Hence we did not further experiment with them.

\bibliographystyle{elsarticle-num}
\bibliography{bibtex}